\documentclass[10pt,twoside,final]{camc}
\usepackage{enumerate}
\usepackage{siunitx}
\DeclareSymbolFont{igrfm}{OML}{cmm}{m}{it}
\def\camcsym#1#2{\DeclareMathSymbol{#1}{\mathord}{igrfm}{"#2}}
\camcsym{\alpha}{0B}\camcsym{\beta}{0C}\camcsym{\gamma}{0D}\camcsym{\delta}{0E}
\camcsym{\epsilon}{0F}\camcsym{\zeta}{10}\camcsym{\eta}{11}\camcsym{\theta}{12}
\camcsym{\iota}{13}\camcsym{\kappa}{14}\camcsym{\lambda}{15}\camcsym{\mu}{16}
\camcsym{\nu}{17}\camcsym{\xi}{18}\camcsym{\pi}{19}\camcsym{\rho}{1A}
\camcsym{\sigma}{1B}\camcsym{\tau}{1C}\camcsym{\upsilon}{1D}\camcsym{\phi}{1E}
\camcsym{\chi}{1F}\camcsym{\psi}{20}\camcsym{\omega}{21}\camcsym{\varepsilon}{22}
\camcsym{\vartheta}{23}\camcsym{\varpi}{24}\let\varrho\rho\let\varsigma\sigma
\camcsym{\varphi}{27}\camcsym{\Gamma}{00}\camcsym{\Delta}{01}\camcsym{\Theta}{02}
\camcsym{\Lambda}{03}\camcsym{\Xi}{04}\camcsym{\Pi}{05}\camcsym{\Sigma}{06}
\camcsym{\Upsilon}{07}\camcsym{\Phi}{08}\camcsym{\Psi}{09}\camcsym{\Omega}{0A}
\let\le\leqslant\let\ge\geqslant

\newcommand{\edmt}[3][\footnotesize]{
\par\noindent{#1\sffamily\bfseries #2~}{#1\rmfamily #3}\vskip.6\baselineskip}

\def\funding#1{\edmt{Funding}{#1}}

\def\complia#1{\edmt[]{Compliance with Ethical Standards}{}\edmt{Conflict of Interest}{#1}}

\AtBeginDocument{%
\def\D{\mathrm{d}}
\newcommand{\ignore}[1]{}
}

\graphicspath{{figs/}}

\begin{document}

\title{A fully discrete Active Flux method for the Euler equations
	comparing different truly multi-dimensional evolution operators}
\author[addressref={aff1}, corref, email={christiane.helzel@hhu.de}]{Christiane~Helzel}%
\author[addressref={aff1},email={amelie.porfetye@hhu.de}]{Amelie~Porfetye}%
\address[id=aff1]{Institute of Mathematics, Heinrich-Heine-University D\"usseldorf, D\"usseldorf, Germany}%
\maketitle
\abstract{%
This article builds upon our recently published fully discrete Cartesian
grid Active Flux method for the Euler equations by introducing a new evolution
operator for the linearized Euler equations.
The evolution of the point value degrees of freedom located at grid
cell boundaries and used to compute numerical fluxes is a key
component of fully discrete Active Flux methods.
Our methods are based on local linearizations of the Euler equations
and the use of truly multi-dimensional evolution operators for 
these linearized problems. 
Here, we propose to solve the linearized Euler equations in moving coordinates,
reducing them to acoustics, which can be solved exactly.
We compare this approach with our previous one, in which
we used approximate evolution operators derived using
the method of bicharacteristics.
The moving-grid approach is closely related to recent methods
proposed by Barsukow as well as Duraisamy. Our choice of local
linearization, as well as correction of linearization errors, enables
us to construct third-order accurate methods for smooth
solutions of the nonlinear Euler equations.
Numerical results illustrate the performance of these methods
for different flow regimes, ranging
from discontinuous solution structures with shock waves to vortex
structures in the low Mach number regime.}


\section{Introduction}

Motivated by the work of P.L.~Roe and several of his former students,
particular Eymann, Fan, Maeng and Samani \cite{article:ER2011a,article:ER2011b,article:ER2013,FR2015,article:RMD2018,samani2023acoustics}, the development 
of Active Flux methods has led to an active and growing
research field. 
Roe \cite{article:Roe2017,article:Roe2021,roe2024musings} proposed to redesign the core components of finite volume 
methods in computational fluid mechanics to obtain fully discrete,
truly multi-dimensional methods with compact stencils that provide
accurate approximations even on relatively coarse grids.

Active Flux methods use cell average values of the conserved
quantities and additional point values along the grid cell boundaries
 as
degrees of freedom.
The point values are used for the
reconstruction as well as the computation of numerical
fluxes. Active Flux methods use globally continuous
piecewise quadratic reconstructions that preserve cell
averages. Such reconstructions have been considered both for triangular grids \cite{article:ER2013}
as well as Cartesian grids \cite{article:BHKR2019,article:HKS2019}.
Here we restrict our considerations to Cartesian grids.
Point values are evolved in time
using truly multi-dimensional evolution operators. They serve as nodes
in quadrature formulas, typically Simpson's rule, for the
numerical fluxes of the finite volume method.

Linear acoustics serves as the prime example to illustrate the benefits of
fully discrete Active Flux methods
\cite{article:ER2013,FR2015,article:BHKR2019,article:Barsukow2020}.
For these methods the point values are evolved using the exact
evolution operator. Barsukow et
al. \cite{article:BHKR2019,article:Barsukow2020} showed that fully
discrete Cartesian grid Active Flux methods for acoustics are vorticity preserving. 
Much earlier, Morton and Roe \cite{article:MR2001} already
pointed out that anomalous numerical
solutions of the Euler equations of gas dynamics, such as the carbuncle
phenomena or the appearance of unphysical Mach stems, appear in
connection with an unphysical vorticity. They conjectured that a
numerical method which controls vorticity is less susceptible to these
unphysical solution structures.
In recent years, research on Active Flux methods has focused on
extending these methods to nonlinear hyperbolic systems, 
particularly the Euler equations of gas dynamics. Eymann and Roe 
proposed splitting methods that approximate transport and acoustics
separately \cite{article:ER2013}. Results are also documented in \cite{PhD:Fan2017}.
Abgrall et al. \cite{article:AB2023,abgrall2025semi} and
Duan et al. \cite{duan2025active}
developed semi-discrete generalized Active Flux methods, which use the
same degrees of freedom and therefore allow compact discretizations of spatial
derivatives. These methods have been shown to be capable of
approximating the Euler
equations in the low Mach number regime. In addition, the recent development
of bound-preserving limiting concepts enabled the application of
these methods for challenging test problems with strong shock waves.

Our own group focuses on developing fully discrete Active
Flux methods and exploring potential evolution operators for point
value updates. In \cite{article:HKS2019}, we proposed an ADER approach for the
evolution of the point values in methods which use the same degrees
of freedom as Active Flux methods. For linear
one-dimensional problems this method is equivalent with Active Flux
methods that use exact evolution
operators. For two-dimensional problems the ADER approach does not use all the
directions of wave propagation. 
In \cite{article:CHL2024}, we proposed Active Flux
methods for linear hyperbolic systems, in particular acoustics and the
linearized Euler equations, 
which use a multi-dimensional evolution operator derived via the method of
bicharacteristics in earlier work by Luk\'a\v{c}ov\'a et al.
\cite{article:LMW2000,article:LSW2002}. 
The evolution operators for
the point value update are third-order accurate and in connection with
the Active Flux method lead to third-order accurate approximations
which are stable for time steps corresponding to CFL$\le
0.279$. Previously, see \cite{article:PTCHL2025},
new third-order accurate evolution operators have
been derived which improve the stability and allow time steps up
to CFL$\le0.44$. For acoustics, we compared the resulting fully
discrete Active Flux methods using either the exact evolution operator
or the third-order accurate operators proposed by the method of
bicharacteristics.
While for most problems the accuracy is comparable, the approximate evolution operators derived via the method of
bicharacteristics do not exactly preserve vorticity.
In \cite{article:CHP2025}, we presented a fully discrete Active Flux
method for the nonlinear Euler equations, which is based on local linearizations
and the use of the method of bicharacteristics for the evolution of
point values. We achieved third-order accuracy for smooth problems by
carefully choosing the local linearization, which allowed us to
correct the linearization error. Furthermore, a bound-preserving
limiting of point and cell average values was introduced, enabling
the method to be applied to problems with discontinuous solution structures.
Fully discrete Active Flux methods are well suited for adaptively refined grids. For Cartesian grids, such an approach was explored in \cite{article:CCH2023} using the ForestClaw software \cite{misc:ForestClaw}.

Recently, Barsukow \cite{article:Barsukow2025} proposed an additive splitting
approach for the evolution of the point values in fully discrete
Active Flux methods, solving the
acoustic equations followed by a nonlinear transport. While this
splitting introduces a splitting error which formally reduces the method to
second-order, numerical results on coarse grids show convergence rates
that are comparable with those seen by third-order accurate
methods. Furthermore, numerical results indicate good stability
properties and excellent performance, in particular when applied in the low
Mach number regime.
While his approach aims to approximate a nonlinear evolution for each
point value, it motivated us to evolve point values using a moving
grid approach for the linearized Euler equations. In combination with
the exact evolution operator for acoustics, this leads to a new exact
solver for the linearized Euler equations. Using the local linearization and the correction of the linearization error as proposed in \cite{article:CHP2025}, we obtain a new third-order accurate method for the nonlinear Euler equations. We presented the moving-grid approach at the W\"urzburg conference in March 2026. In
the meantime, we became aware of a new preprint by Duraisamy
\cite{duraisamy2026}
that also uses a moving coordinate approach. 

In this paper, we review our fully discrete
Active Flux method for the Euler equations, using either the truly
multi-dimensional evolution operator derived via the method of
bicharacteristics or the new moving-grid approach.
We use the same local linearization and correction of the
linearization error for both versions. For problems with shock waves
we use the limiting strategy described in \cite{article:CHP2025}. 

The rest of the paper is organized as follows: In Section \ref{sec:2}, we
provide a brief review of the main components of
the fully discrete Active Flux method from \cite{article:CHP2025}. In
Section \ref{sec:3}, we present the two different multi-dimensional
evolution operators for
the evolution of point values for the linearized Euler
equations. In Section \ref{sec:4}, 
we discuss the performance of the resulting Active Flux methods for
several linear test problems. Finally, in Section \ref{sec:5}, we present the
results obtained using these different evolution operators in a fully discrete
Active Flux method for the nonlinear Euler equations. 

\section{Review of the fully discrete Active Flux method for
	Euler}
\label{sec:2}

We give a brief review of the fully discrete Active flux method,
a more detailed description can be found in \cite{article:CHP2025}.
We consider the two-dimensional Euler equations of gas dynamics in conservative form
\begin{equation}
	\partial_t \mathbf{q} + \partial_x \mathbf{f}(\mathbf{q}) + \partial_y \mathbf{g}(\mathbf{q}) = \mathbf{0},
\end{equation}
with 
$\mathbf{q}  = (\rho, \rho u, \rho v, E)^T$, 
$\mathbf{f}(\mathbf{q})  = (\rho u, \rho u^2 + p, \rho u v, u(E+p))^T$,
$\mathbf{g}(\mathbf{q})  = (\rho v, \rho u v,\\
 \rho v^2 + p, v(E+p))^T$
and ideal gas equation of state 
$$
E =\frac{p}{\gamma -1} + \frac{1}{2} \rho (u^2 + v^2).
$$
Here, $\rho$ is the fluid density, $p$ is the pressure, $u$ and $v$
are the velocity components in the $x$- and $y$-directions, $E$
is the total energy per unit volume, and $\gamma$ is the ratio of
specific heats.

On a Cartesian grid, the cell average values of the conserved
quantities are evolved in time using the finite volume approach
\begin{equation}\label{eqn:fvm}
	\bar{Q}_{i,j}^{n+1} = \bar{Q}_{i,j}^n - \frac{\Delta t}{\Delta x} \left(
	F_{i+\frac{1}{2},j} - F_{i-\frac{1}{2},j}\right) - \frac{\Delta
		t}{\Delta y} \left(G_{i,j+\frac{1}{2}}-G_{i,j-\frac{1}{2}}\right),
\end{equation}
where the fluxes are approximated using Simpson's rule, i.e.,
\begin{equation} \label{eqn:fluxF}
	\begin{split}
		& \frac{1}{\Delta t \Delta y}
		\int_{t_n}^{t_{n+1}}
		\int_{y_{j-\frac{1}{2}}}^{y_{j+\frac{1}{2}}}
		\mathbf{f}(\mathbf{q}(x_{i+\frac{1}{2}},y,t))
		\,{\rm d}y\,{\rm d}t \\
		& \approx \frac{1}{36}
		\Bigl(
		\mathbf{f}(Q_{i+\frac{1}{2},j-\frac{1}{2}}^n)
		+4\mathbf{f}(Q_{i+\frac{1}{2},j}^n)
		+\mathbf{f}(Q_{i+\frac{1}{2},j+\frac{1}{2}}^n)
		\\
		& \qquad
		+4\mathbf{f}(Q_{i+\frac{1}{2},j-\frac{1}{2}}^{n+\frac{1}{2}})
		+16\mathbf{f}(Q_{i+\frac{1}{2},j}^{n+\frac{1}{2}})
		+4\mathbf{f}(Q_{i+\frac{1}{2},j+\frac{1}{2}}^{n+\frac{1}{2}})
		\\
		& \qquad
		+\mathbf{f}(Q_{i+\frac{1}{2},j-\frac{1}{2}}^{n+1})
		+4\mathbf{f}(Q_{i+\frac{1}{2},j}^{n+1})
		+\mathbf{f}(Q_{i+\frac{1}{2},j+\frac{1}{2}}^{n+1})
		\Bigr)
		\\
		& =: F_{i+\frac{1}{2},j}
	\end{split}
\end{equation}
and analogously for $G_{i,j+\frac{1}{2}}$.

In Active Flux methods the degrees of freedom are cell average values
of the conserved variables, denoted by $\bar{Q}_{i,j}^n$,
and point values along the grid cell boundary, denoted by
$Q_{i-\frac{1}{2},j-\frac{1}{2}}^n$, $Q_{i,j-\frac{1}{2}}^n$,
$Q_{i+\frac{1}{2},j-\frac{1}{2}}^n$, $Q_{i+\frac{1}{2},j}^n$,
$Q_{i+\frac{1}{2},j+\frac{1}{2}}^n$, $Q_{i,j+\frac{1}{2}}^n$,
$Q_{i-\frac{1}{2},j+\frac{1}{2}}^n$ and $Q_{i-\frac{1}{2},j}^n$. We
denote with $U$ point and cell average values in primitive variables. Point values in primitive variables can straight forwardly
be computed from point values in conservative variables. Furthermore, cell average values can be transformed between conservative and primitive variables with the required high-order accuracy by also making use of the point values.

A globally continuous, piecewise quadratic reconstruction
$\mathbf{u}_{rec}^n(x,y)$ is computed
from point and cell average values in primitive variables at time
$t_n$, using basis functions presented, for example, in
\cite{article:BHKR2019,article:HKS2019}.
The crucial step is the computation of point values at time
$t_{n+\frac{1}{2}}$ and $t_{n+1}$ that are needed for the numerical
fluxes (\ref{eqn:fluxF}) as well as the reconstruction at the next time step.

\subsection{Evolution of the point values using local linearization}\label{sec:2-1}
The evolution of the point values is based on the
quasilinear form of the equations in primitive variables, i.e., 
\begin{equation}\label{eqn:quasilin}
	\partial_t \mathbf{u} + A(\mathbf{u}) \partial_x \mathbf{u} + B(\mathbf{u}) \partial_y \mathbf{u} = \mathbf{0}, 
\end{equation}
with the vector of primitive variables $\mathbf{u} = (\rho, u, v, p)^T$ and
\begin{equation}\label{eon:quasilinAB}
	A(\mathbf{u}) := \left( \begin{array}{cccc}u & \rho & 0 & 0 \\ 0 & u &
		0 & 1/\rho \\ 0 & 0 & u & 0 \\ 0 & \gamma p & 0 & u\end{array}\right), \quad  B(\mathbf{u}) :=  \left( \begin{array}{cccc}v & 0 & \rho & 0 \\ 0 & v & 0 & 0 \\ 0 & 0 & v & 1/\rho \\ 0 & 0 & \gamma p & v \end{array}\right).
\end{equation}
We seek third-order accurate approximations of (\ref{eqn:quasilin})
at the positions of the point values and at the
times $t_{n+\frac{1}{2}}=t_{n} +\Delta t/2$ and $t_{n+1} = t_n+\Delta
t$, using the initial condition $\mathbf{u}(x,y,t_n) = \mathbf{u}_{rec}^{n}(x,y)$. 

Approximate solutions of the nonlinear problem (\ref{eqn:quasilin})
are obtained 
by performing, at each point value, a local linearization and solving
a linear hyperbolic system

\begin{equation}\label{eqn:linEuler}
	\partial_t \mathbf{v} + A \partial_x \mathbf{v} + B \partial_y \mathbf{v} = \mathbf{0},
\end{equation}
with initial values $\mathbf{v}(x,y,t_n) = \mathbf{u}_{rec}^{n}(x,y)$.
Here, $\mathbf{v}$ is again the vector of primitive variables and
$A:=A(\mathbf{u}')$ and $B:=B(\mathbf{u}')$ are the matrices from
(\ref{eon:quasilinAB}) evaluated at a constant state
$\mathbf{u}'$ used for the local linearization. In Section \ref{sec:3}, we will introduce two truly
multi-dimensional evolution operators for linear hyperbolic
systems. One of these, which was also used in
\cite{article:CHL2024,article:CHP2025}, provides an 
approximation of solutions to the linear hyperbolic system that is
third-order accurate. The second
method provides an exact solution to the linearized Euler
equations. Note that, in both cases, the piecewise quadratic
reconstruction limits
the accuracy of the method for linear hyperbolic problems to third-order.

A linearization error arises when the solution to
equation (\ref{eqn:quasilin}) is replaced by the solution to equation
(\ref{eqn:linEuler}). The precise form of this error depends on the
choice of $\mathbf{u}'$, but no linearization provides
third-order accurate results. \ignore{In \cite[Theorem 2.1]{article:CHP2025}
we provide the leading term of the truncation error of
$\mathbf{u}(\bar{x},\bar{y},t_n+\tau) -
\mathbf{v}(\bar{x},\bar{y},t_n+\tau)$ when a linearization
around $\mathbf{u}'=\mathbf{u}(\bar{x},\bar{y},t_n+\tau/2)$ is used to
define the linearized Euler equations. }
In \cite[Theorem 2.1]{article:CHP2025}, we derived the leading term of the linearization error
\begin{align*}
\mathbf{u}(\bar{x},\bar{y},t_n+\tau) &-
\mathbf{v}(\bar{x},\bar{y},t_n+\tau) = \\
     \frac{1}{2}  &\tau^2\left( A(\mathbf{u}) \frac{\partial A(\mathbf{u})}{\partial
    \mathbf{u}} \cdot (\partial_x \mathbf{u},\partial_x \mathbf{u}) +
  A(\mathbf{u}) \frac{\partial B(\mathbf{u})}{\partial
    \mathbf{u}} \cdot (\partial_x \mathbf{u},\partial_y \mathbf{u}) \right.\\
\qquad \quad +& \left.  B(\mathbf{u}) \frac{\partial A(\mathbf{u})}{\partial
    \mathbf{u}} \cdot (\partial_y \mathbf{u},\partial_x \mathbf{u}) +
  B(\mathbf{u}) \frac{\partial B(\mathbf{u})}{\partial
    \mathbf{u}} \cdot (\partial_y \mathbf{u},\partial_y \mathbf{u}) \right) \bigg\vert_{(\bar{x},\bar{y},t_n)}\\
    +& \mathcal{O}(\tau ^3),
\end{align*}
assuming that
$\mathbf{u}'=\mathbf{u}(\bar{x},\bar{y},t_n+\tau/2)$ is used for the local linearization.
A third-order accurate
approximation of the point value
$\mathbf{u}(\bar{x},\bar{y},t_n+\tau)$ can be obtained by adding a
correction term, which discretizes the leading term of the truncation
error as explained in \cite[Corollary 2.2]{article:CHP2025}.

\ignore{
When approximating discontinuous solution structures, limiting
may be necessary to preserve positivity of density and pressure, and to
avoid unphysical oscillations. In \cite{article:CHP2025}, we presented
a limiting of point as well as cell average values, blending the
third-order Active Flux approximation discussed above with a
low-order, bound-preserving method. Our limiting strategy is closely
related to the limiting concepts developed for semi-discrete Active
Flux methods, see for example Duan et al. \cite{duan2025active} and
Abgrall et al. \cite{abgrall2026bound}. A key difference from
semi-discrete methods is that our fully discrete method also requires
the state used for
the local linearization to be limited. See \cite[Section 3]{article:CHP2025} for a detailed
description of our limiting procedure.}

\subsection{Limiting }\label{sec:limiting} 
When approximating discontinuous solution structures, limiting
may be necessary to preserve positivity of density and pressure, and to
avoid unphysical oscillations. 
Our limiting strategy is 
related to the limiting concepts developed for semi-discrete Active
Flux methods, see for example Duan et al. \cite{duan2025active} and
Abgrall et al. \cite{abgrall2026bound}. A key difference from
semi-discrete methods is that our fully discrete method 
focuses on limiting point values rather than fluxes.
See \cite[Section 3]{article:CHP2025} for a detailed
description of our limiting procedure. 

In the following, we briefly summarize the main components of the limiting strategy.
An important component is a shock indicator of the form
$$\theta= \exp(-\kappa \phi^{(1)}\phi^{(2)}),$$ which detects shocks and also reflects their strength.
$\phi^{(1)}$ is based on normalized finite-difference approximations of second derivatives of the pressure and therefore detects strong pressure variations while $\phi^{(2)}$ is based on finite-difference approximations of the divergence and vorticity of the velocity field and $\kappa > 0$ is a fixed parameter. For all of our numerical experiments, we choose $\kappa = 2$. Consequently, $\theta$ is close to one in smooth regions and decreases towards zero in the vicinity of strong shocks.

 We observed that limiting the state $\mathbf{u}'$, which is used to compute the matrices $A$ and $B$ of the linearized Euler equations, in the vicinity of shocks improves the method's accuracy.
 A suitable state for the linearization is obtained as a blend of $\mathbf{u}'$, as used in the third-order method, and neighboring cell averages from the previous time step. We denote by $\theta_{i-\frac{1}{2},j-\frac{1}{2}}$ the shock indicator associated with the grid point $x_{i-\frac{1}{2},j-\frac{1}{2}}$. The corresponding state used for the local linearization is then defined by
\begin{equation}\label{eq:limiting_linearisation}
\begin{split}
 U'_{i-\frac{1}{2},j-\frac{1}{2}}
&= s\left(\theta_{i-\frac{1}{2},j-\frac{1}{2}}\right)
U_{i-\frac{1}{2},j-\frac{1}{2}}(t_n+\frac{\tau}{2})\\
&\quad + \left(1-s\left(\theta_{i-\frac{1}{2},j-\frac{1}{2}}\right)\right)
\frac{1}{4}\left(\bar{U}_{i-1,j-1}^n+\bar{U}_{i,j-1}^n+\bar{U}_{i-1,j}^n+\bar{U}_{i,j}^n \right), 
\end{split}
\end{equation}
where $s:[0,1]\rightarrow[0,1]$ is a blending function.
Here, $U_{i-\frac{1}{2},j-\frac{1}{2}}(t_n+\frac{\tau}{2})$ is the point value in primitive variables at an intermediate time level as suggested by our accuracy study provided in \cite[Theorem 2.1]{article:CHP2025} and outlined above. 

We observed that the new evolution operator introduced in Section \ref{sec:moving} requires slightly stronger limiting of the state vector used for the local linearization than the EG2 evolution operators employed in \cite{article:CHP2025}. 
 Therefore, rather than relying solely on the 
 shock indicator associated with a given grid point, 
 we also take into account the values at the neighboring corners. More precisely, in \eqref{eq:limiting_linearisation} we replace $\theta$ by a modified indicator $\theta'$, defined by $\theta_{i-\frac{1}{2},j-\frac{1}{2}}' := \theta_{i-\frac{1}{2},j-\frac{1}{2}} $, $\theta_{i,j-\frac{1}{2}}' := \min \left\{\theta_{i-\frac{1}{2},j-\frac{1}{2}},\theta_{i+\frac{1}{2},j-\frac{1}{2}} \right\}$ and $\theta_{i-\frac{1}{2},j}' := \min \left\{\theta_{i-\frac{1}{2},j-\frac{1}{2}},\theta_{i-\frac{1}{2},j+\frac{1}{2}} \right\}$. For consistency, this form of the shock indicator is used in all computations presented in this paper that involve limiting.

 The shock indicator is also employed in the limiting of the point values as explained in detail in \cite[Section 3]{article:CHP2025}. First, point values that are not bound-preserving are replaced by low-order but bound-preserving approximations. In a second step, we compute a blend of high- and low-order point values, i.e.,
 $$ U^{(lim)}  := \theta U^{(ho)} + \left(1-\theta\right)U^{(lo)},$$
 where the amount of blending is again determined by the shock indicator.
 This additional limiting reduces unphysical oscillations near discontinuities. 
 
 The final step of our limiting strategy is the limiting of the numerical fluxes to obtain bound-preserving cell averages. The general idea, in analogy to the limiting of point values, is to replace the fluxes
in the finite volume method \eqref{eqn:fvm} with a convex linear combination of the high-order
flux and a local Lax-Friedrichs flux. For the numerical experiments presented in this work, however, no flux limiting was needed.
\section{Truly multi-dimensional evolution operators
	for the linearized Euler equations}\label{sec:3}
In this Section, we present two different evolution operators for the linearized
Euler equations.
Both will be compared in numerical simulations for the linearized as
well as nonlinear Euler equations. The first, briefly reviewed in
Section \ref{sec:EG2} is a third-order accurate
operator proposed by Lukacova et al. \cite{article:LSW2002}
and used in fully discrete Active Flux methods \cite{article:CHL2024,article:CHP2025}.
The second operator, described in Section \ref{sec:moving},
uses a moving-grid approach and describes an
exact evolution of point values for the linearized Euler equations.

We consider linear hyperbolic systems of the form
\begin{equation}\label{eqn:linEuler2}
	\partial_t \mathbf{v} + A \partial_x \mathbf{v} + B \partial_y
	\mathbf{v} = \mathbf{0},
\end{equation}
with $\mathbf{v}:\mathbb{R}^2 \times \mathbb{R}^+ \rightarrow
\mathbb{R}^4$, where the components 
$\mathbf{v} = (\rho, u, v, p)^T$ represent density,
velocity in $x$- and $y$- direction,
and pressure. The constant matrices have the form $A:=A(\mathbf{u}')$
and $B:=B(\mathbf{u}')$ with $A$ and $B$ as described in equation
(\ref{eon:quasilinAB}).

This system is obtained by
linearizing the compressible Euler equations about a constant
background state $({\rho}',{u}',{v}',{p}') \in
\mathbb{R}^4$, with ${\rho}',{p}' >0$, 
and the speed of sound 
$c':=\sqrt{\gamma p'/\rho'}$.
Note that in this Section, we are assuming a constant background state, denoted by $\mathbf{u}'$, to linearize the equations across the entire domain. In order to be consistent with the notation used in Section \ref{sec:2-1}, we denote the conservative variable of the linear system with $\mathbf{v}$. 
\subsection{Evolution of point values using the method of
	bicharacteristics}\label{sec:EG2}
\begin{figure}[htb]
	\centerline{\includegraphics[width=0.7\textwidth]{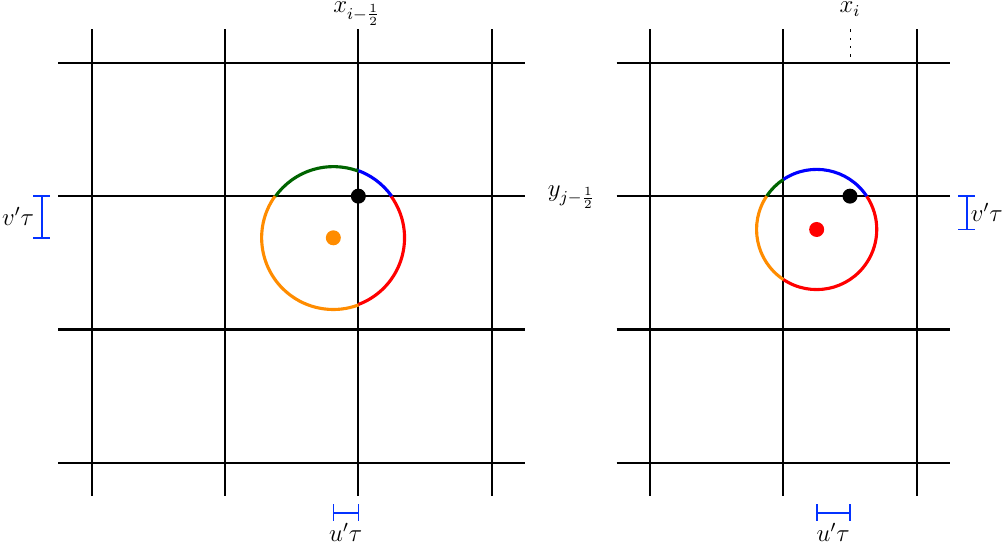}}
	\caption{\label{fig:EG2}Schematic description of the EG2 evolution operator for a
		point value at a grid cell corner and a point value at the
		midpoint of a grid cell interface for $u' > 0$ and $v' > 0$.}
\end{figure}

Third-order accurate evolution operators for the linearized Euler
equations have been derived by Luk\'a\v{c}ov\'a et
al. \cite{article:LSW2002} using the method of bicharacteristics. Let
$P:=(x,y,t_n+\tau)$ be a point at which we want to compute the solution.
Furthermore, we define $P':= (x-{u}'\tau, y-{v}'\tau,t_n)$ and
$Q(\theta):= (x-({u}'-c\cos(\theta))\tau,
y-({v}'-c\sin(\theta))\tau,t_n)$. Then this evolution operator can
be expressed in the form
\begin{equation}\label{eqn:EG2lineuler}
	\begin{alignedat}{2}
		\rho(P) &=
		\rho({P'})-2\frac{p({P'})}{c'^2}
		+\frac{1}{\pi}\int_{0}^{2\pi}
		\frac{p(Q(\theta))}{c'^2}
		-\frac{\rho'}{c'}u(Q(\theta))\cos(\theta)\\
		&\qquad
		-\frac{\rho'}{c'}v(Q(\theta))\sin(\theta)
		\,\D\theta
		+{\cal O}(\Delta t^3),\\
		u(P) &=
		\frac{1}{\pi}\int_{0}^{2\pi}
		-\frac{p(Q(\theta))}{\rho'c'}\cos(\theta)
		+u(Q(\theta))
		\left(2\cos^2(\theta)-\frac{1}{2}\right)\\
		&\qquad
		+2v(Q(\theta))\sin(\theta)\cos(\theta)
		\,\D\theta
		+{\cal O}(\Delta t^3),\\
		v(P) &=
		\frac{1}{\pi}\int_{0}^{2\pi}
		-\frac{p(Q(\theta))}{\rho'c'}\sin(\theta)
		+2u(Q(\theta))\sin(\theta)\cos(\theta)\\
		&\qquad
		+v(Q(\theta))
		\left(2\sin^2(\theta)-\frac{1}{2}\right)
		\,\D\theta
		+{\cal O}(\Delta t^3),\\
		p(P) &=
		-p({P'})
		+\frac{1}{\pi}\int_{0}^{2\pi}
		p(Q(\theta))
		-\rho'c'u(Q(\theta))\cos(\theta)\\
		&\qquad
		-\rho'c'v(Q(\theta))\sin(\theta)
		\,\D\theta
		+{\cal O}(\Delta t^3).
	\end{alignedat}
\end{equation}
Active Flux methods for linear hyperbolic systems using these kinds of
evolution operators have been proposed in \cite{article:CHL2024}. An
eigenvalue analysis for acoustics revealed stability for CFL$\le 0.279$.
We observed that the same time step restriction is needed for the
linearized Euler equations. Figure \ref{fig:EG2} provides a schematic illustration of the EG2 evolution operator for the update of a corner point value and an edge midpoint.
 
\subsection{Evolution of point values using moving-grid approach}
\label{sec:moving}
\begin{figure}
	\centerline{\includegraphics[width=0.7\textwidth]{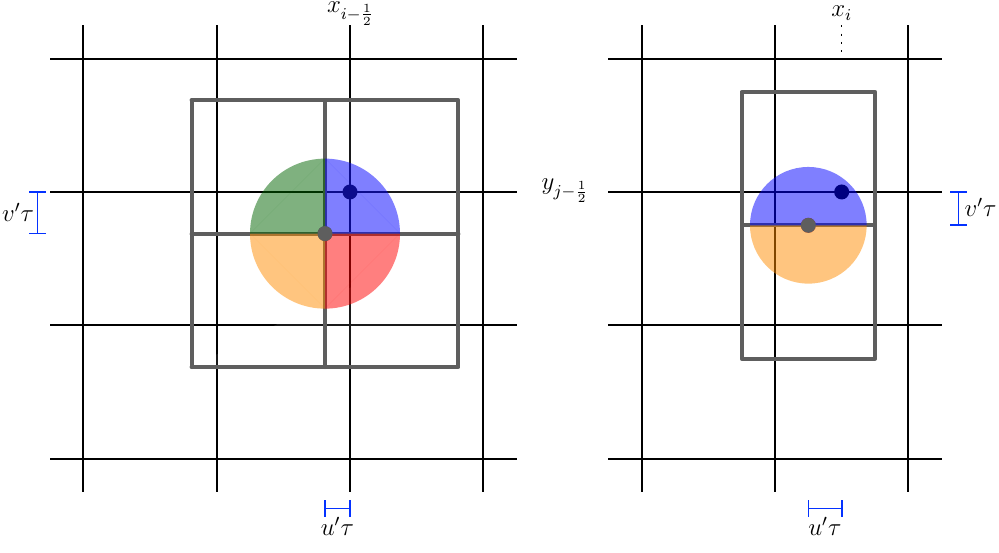}}
	\caption{\label{fig:moving}Illustration of the moving-grid approach for the update of a
		corner point value and the midpoint of the bottom interface for $u' > 0$ and $v' > 0$.} 
\end{figure}

In our moving-grid approach we consider coordinates $\xi =
x-{u}'t$, $\eta = y-{v}'t$, $\tau = t$. Then the partial
derivatives transform according to $\partial_t = \partial_\tau -
{u}' \partial_\xi - {v}' \partial_\eta$, $\partial_x =
\partial_\xi$, $\partial_y = \partial_\eta$ and the linearized Euler
equations can be written in the form
\begin{equation}\label{eqn:moving-1}
	\begin{split}
		\partial_\tau \rho + {\rho}'\left( \partial_\xi u + \partial_\eta v
		\right) & = 0, \\
		\partial_\tau u + \frac{1}{{\rho}'} \partial_\xi p & = 0, \\
		\partial_\tau v + \frac{1}{{\rho}'} \partial_\eta p & = 0, \\
		\partial_\tau p + {\rho}' c'^2 \left( \partial_\xi u + \partial_\eta
		v \right) & = 0. 
	\end{split}
\end{equation}
Introducing the scaled pressure variable $\tilde{p} =
\frac{p}{{\rho}' c'}$, i.e., $p = {\rho}' c' \tilde{p}$ allows us to rewrite the second, third, and fourth equations of (\ref{eqn:moving-1}) as
\begin{equation}\label{eqn:moving-2}
	\begin{split}
		\partial_\tau u + c \partial_\xi \tilde{p} & = 0,\\
		\partial_\tau v + c \partial_\eta \tilde{p} & = 0, \\
		\partial_\tau \tilde{p} + c' \left( \partial_\xi u + \partial_\eta v
		\right) & = 0.
	\end{split}
\end{equation}
Thus, in these moving coordinates, the variables $\tilde{p}$, $u$, and $v$
form a closed acoustic system, given by (\ref{eqn:moving-2}), that is
independent of the density $\rho$.
However, the rescaled density $\tilde{\rho} =
\frac{c'}{{\rho}'} \rho$ satisfies an evolution equation of the form
\begin{equation}
	\partial_\tau \tilde{\rho} + c' \left( \partial_\xi u + \partial_\eta v
	\right) = 0,
\end{equation}
which coincides with the evolution equation for $\tilde{p}$.
Thus, the evolution of $\tilde{\rho}$ can
be computed from increments of $\tilde{p}$. 
The resulting moving-grid evolution strategy is illustrated in Figure \ref{fig:moving} for two different point values, one located at the corner of a grid cell and one located at the midpoint of a horizontal edge.

A point value located at the midpoint of a vertical edge can be
handled analogously.
Depending on the location of the point value, we compute a new grid
patch consisting of either two or four Cartesian grid cells.
We assign
point and cell average values of $\mathbf{v}$ to this grid patch, which
allows us to compute continuous, piecewise quadratic
reconstructions. Finally, we apply the exact evolution operator for
acoustics to the shifted grid cell patch.

We will now present the details of the point value update for a
corner point located at $(x_{i-\frac{1}{2}},y_{j-\frac{1}{2}})$.
Let $V_{i-\frac{1}{2},j-\frac{1}{2}}^n$ denote the point
values of $\mathbf{v}=(\rho,u,v,p)$ at time $t_n$ on the original
grid. Our goal is to compute
$V_{i-\frac{1}{2},j-\frac{1}{2}}^{n+\frac{1}{2}}$ and $V_{i-\frac{1}{2},j-\frac{1}{2}}^{n+1}$.
Furthermore, let $\mathbf{v}_{rec}^n:\mathbb{R}^2\rightarrow
\mathbb{R}^4$ denote the globally continuous, piecewise quadratic
reconstruction based on cell average values and point values of the original grid
at time $t_n$.

In the first step, we compute the point
$(\bar{x},\bar{y}) = (x_{i-\frac{1}{2}}-u'
\tau,y_{j-\frac{1}{2}}-v'\tau)$
for $\tau \in \{\frac{\Delta t}{2},\Delta t\}$, where
	$(\bar{x},\bar{y})$ corresponds to the point $P'$ introduced in Section~\ref{sec:EG2}, and define the local domain
$$\Omega_{i-\frac{1}{2},j-\frac{1}{2}}:=[\bar{x}-\Delta
x, \bar{x}+\Delta x] \times [\bar{y}-\Delta y, \bar{y}+\Delta y].$$
The domain, $\Omega_{i-\frac{1}{2},j-\frac{1}{2}}$ is discretized by four cells, 
as indicated in Figure \ref{fig:moving} (left).
Next, using
the underlying piecewise quadratic reconstruction
$\mathbf{v}_{rec}^n$ of the original grid,
we evaluate the physical quantities at the corner points
and edge midpoints of the cells in $\Omega_{i-\frac{1}{2},j-\frac{1}{2}}$ and compute the cell average values for its four
cells. In particular, the physical quantities associated with the corner point $(\bar{x},\bar{y})$ are given by $\hat{V}_{i-\frac{1}{2},j-\frac{1}{2}} := \mathbf{v}_{rec}^n(\bar{x},\bar{y})$.
Based on these newly computed point values and cell averages, we construct a continuous piecewise quadratic Active Flux reconstruction on $\Omega_{i-\frac{1}{2},j-\frac{1}{2}}$, denoted by $\hat{\mathbf{v}}_{rec}$.

We now consider the
pressure and velocity components of
$\hat{\mathbf{v}}_{rec}-\hat{V}_{i-\frac{1}{2},j-\frac{1}{2}}$, 
scale the pressure perturbation, and introduce the notation
\begin{equation*}
	(\tilde{p}_0,{u}_0,{v}_0):=\left((\hat{p}_{rec}-\hat{p}_{i-\frac{1}{2},j-\frac{1}{2}})/(\rho'
	c'), \hat{u}_{rec}-\hat{u}_{i-\frac{1}{2},j-\frac{1}{2}},\hat{v}_{rec}-\hat{v}_{i-\frac{1}{2},j-\frac{1}{2}}\right).
\end{equation*}
We then compute  the solution of the 
acoustic system
(\ref{eqn:moving-2}) with piecewise quadratic initial data
$(\tilde{p}_0,{u}_0,{v}_0)$ at the point
$(\bar{x},\bar{y})$ for $\tau \in \{\frac{\Delta t}{2}, \Delta t\}$.
The components of this solution, which can be interpreted as increments with respect to a reference state at $(\bar{x},\bar{y})$, are denoted by
$(\delta \tilde{p}(\tau), \delta u(\tau), \delta v(\tau))$.
We use the exact acoustic solver,
proposed by Barsukow et al. \cite{article:BHKR2019}, which has the form
\begin{equation}\label{eqn:Barsukow}
	\begin{aligned}
		p({\bf x},t) & = \partial_r \left( r M_r \{p({\bf x},0) \} \right) \vert_{r=c t} -
		\frac{1}{c t} \partial_r \left( r^2 M_r \{ {\bf n} \cdot
		\mathbf{u}({\bf x},0)  \} \right) \vert_{r=c t} \\
		\mathbf{u}({\bf x},t) & = \mathbf{u}({\bf x},0) - \frac{1}{c t} \partial_r
		\left( r^2 M_r \{  \mathbf{n} p({\bf x},0) \} \right) \vert_{r = ct},
		\\
		& \qquad + \int_0^{c
			t} \frac{1}{r} \partial_r \left( \frac{1}{r} \partial_r
		\left( r^3 M_r \{ ({\bf n} \cdot {\bf u}({\bf x},0) ) {\bf n}
		\right) - r M_r \{ {\bf u}({\bf x},0) \} \right) \, \D r,
	\end{aligned}
\end{equation}
where $M_r\{f\}$ is the spherical mean of a function $f$ over a disc with radius $r$.
Finally, the components of $V_{i-\frac{1}{2},j-\frac{1}{2}}^{n+1}$ are obtained by adding the evolved perturbations to the reference state $\hat{V}_{i-\frac{1}{2},j-\frac{1}{2}}$. While the velocity perturbations are added directly, the scaled pressure perturbation $\delta\tilde{p}$ is transformed back into the corresponding density and pressure perturbations. Thus, we obtain
\begin{equation}
	\begin{split}
		\rho_{i-\frac{1}{2},j-\frac{1}{2}}^{n+1} & =
		\hat{\rho}_{i-\frac{1}{2},j-\frac{1}{2}} + \delta \tilde{p}(\Delta t) \rho'/c',\\  
		u_{i-\frac{1}{2},j-\frac{1}{2}}^{n+1} & =
		\hat{u}_{i-\frac{1}{2},j-\frac{1}{2}} + \delta u(\Delta t), \\
		v_{i-\frac{1}{2},j-\frac{1}{2}}^{n+1} & =
		\hat{v}_{i-\frac{1}{2},j-\frac{1}{2}} + \delta v(\Delta t), \\
		p_{i-\frac{1}{2},j-\frac{1}{2}}^{n+1} & =
		\hat{p}_{i-\frac{1}{2},j-\frac{1}{2}} + \delta \tilde{p}(\Delta t) 
		\rho' c'. \\
	\end{split}
\end{equation}
The components of $V_{i-\frac{1}{2},j-\frac{1}{2}}^{n+\frac{1}{2}}$
are computed analogously.
For the evolution of the point values at the edge midpoints,
i.e., $V_{i-\frac{1}{2},j}$ and $V_{i,j-\frac{1}{2}}$, we use a moving
grid approach where the local domain consists of two grid cells as indicated in Figure
\ref{fig:moving} (right).
Note that the accuracy of the moving-grid approach for the
linearized Euler equations is limited solely by the accuracy of the
Active Flux reconstruction, which limits the accuracy of the point values to third-order. 

The Cartesian grid Active Flux method for acoustics with the exact evolution operator is stable for time steps satisfying CFL$\le0.5$. We used the same time step restriction for all of our simulations performed for the linearized as well as nonlinear Euler equations. Thus, the moving-grid approach allows larger time steps than the fully discrete Active Flux method with the EG2 evolution operator. 

\begin{remark}
	A slightly simpler implementation might construct a single grid cell
	with the shifted point in the center. Our choice of the implementation
	allowed us to make direct use of our previous implementation of the
	exact acoustics solver on a grid that has the same resolution as the
	original grid.
	We plan to explore alternative implementations and their
	impact on stability and accuracy in the future.
\end{remark}

\section{Numerical results for the linearized Euler equations}
\label{sec:4}
In this Section, we present numerical results of the fully discrete
Active Flux method for the linearized Euler equations, using either the
approximate evolution operator from Section \ref{sec:EG2} or the
exact evolution operator from Section \ref{sec:moving}.
In the linear case, we obtain a finite volume method of the form given
in equation 
(\ref{eqn:fvm}) by replacing the terms $\mathbf{f}(\mathbf{q})$  in
(\ref{eqn:fluxF}) by $A \mathbf{q}$ and, analogously,
$\mathbf{g}(\mathbf{q})$ by $B \mathbf{q}$. See also
\cite{article:CHL2024} for a detailed description of the Active Flux
method in the linear case. 

We begin with a test case that is periodic in time, allowing a convergence study to be carried out by simply comparing the solution with the initial data. 
\begin{example}\label{ex:linEu-convergence}
	We consider the linearized Euler equations with $\gamma = 1.4$,
	$\rho'=1$, $u' = v' = 1$, $p'=1/\gamma$, and initial values of
	the form
	\begin{equation}
		\begin{split}
			\rho(x,y,0) & = \cos(2 \pi x), \\
			u(x,y,0) & = -\sin(2 \pi x) + \sin(2\pi y), \\
			v(x,y,0) & = \sin(2 \pi x) + \sin(2 \pi y), \\
			p(x,y,0) & = \cos(2 \pi x),
		\end{split}
	\end{equation}
	on the domain $[0,1]\times[0,1]$ using periodic boundary conditions in
	the $x$- and $y$- direction. At time $t=1$ the exact solution agrees
	with the initial conditions.
\end{example}
In Tables \ref{table:linEu-con-moving} and \ref{table:linEu-con-EG2},
we show the results of a numerical convergence study at $t=1$. Both methods
converge with third-order and the size of the error is
comparable. The Active Flux method with
exact evolution operator allows to use larger time steps.
\begin{table}[htb!]
	\caption{\label{table:linEu-con-moving}Numerical convergence study for Example
		\ref{ex:linEu-convergence} using the fully discrete Active Flux
		method with the exact evolution operator described in Section
		\ref{sec:moving}. The time step restriction corresponds to
 		CFL$\le 0.5$.}
	
	
	\sisetup{
		scientific-notation = true,
		round-mode = places
	}
	
	\begin{center}
		\begin{minipage}{0.9\textwidth}
			
			\begin{tabular}{
					*1{S[
						table-column-width=0.5cm,
						table-text-alignment=center
						]}
					*4{S[
						table-column-width=1.2cm,
						table-text-alignment=center
						]}
					*4{S[
						table-column-width=0.5cm,
						table-text-alignment=center,
						round-precision=4
						]}
				}
				
				\toprule
				
				{Res.}
				& \multicolumn{4}{c}{Error in}
				& \multicolumn{4}{c}{EOC in} \\
				
				\midrule
				
				& {$\bar{\rho}$}
				& {$\bar{u}$}
				& {$\bar{v}$}
				& {$\bar{p}$}
				& {$\bar{\rho}$}
				& {$\bar{u}$}
				& {$\bar{v}$}
				& {$\bar{p}$} \\
				
				\midrule

				{64}
				& \num{4.5125772384423883e-04}
				& \num{4.4043628612257973e-04}
				& \num{3.4243117054858028e-04}
				& \num{4.5125772384423883e-04}
				& {---}
				& {---}
				& {---}
				& {---} \\
				
				{128}
				& \num{5.6617484476779194e-05}
				& \num{5.5265724405451904e-05}
				& \num{4.2976773335026655e-05}
				& \num{5.6617484476779194e-05}
				& {2.99}
				& {2.99}
				& {2.99}
				& {2.99} \\
				
				{256}
				& \num{7.0854758160756855e-06}
				& \num{6.9164685064347308e-06}
				& \num{5.3787764875558739e-06}
				& \num{7.0854758160756855e-06}
				& {3.00}
				& {3.00}
				& {3.00}
				& {3.00} \\
				
				\bottomrule
				
			\end{tabular}
			
		\end{minipage}
	\end{center}
	
\end{table}


\begin{table}[htb!]
	\caption{\label{table:linEu-con-EG2}Numerical convergence study for Example
		\ref{ex:linEu-convergence} using the fully discrete Active Flux
		method with the EG2 evolution operator described in Section
		\ref{sec:EG2}. The time step restriction corresponds to
		CFL$\le 0.279$.}
	
	
	\sisetup{
		scientific-notation = true,
		round-mode = places
	}
	
	\begin{center}
		\begin{minipage}{0.9\textwidth}
			
			\begin{tabular}{
					*1{S[
						table-column-width=0.5cm,
						table-text-alignment=center
						]}
					*4{S[
						table-column-width=1.2cm,
						table-text-alignment=center
						]}
					*4{S[
						table-column-width=0.5cm,
						table-text-alignment=center,
						round-precision=4
						]}
				}
				
				\toprule
				
				{Res.}
				& \multicolumn{4}{c}{Error in}
				& \multicolumn{4}{c}{EOC in} \\
				
				\midrule
				
				& {$\bar{\rho}$}
				& {$\bar{u}$}
				& {$\bar{v}$}
				& {$\bar{p}$}
				& {$\bar{\rho}$}
				& {$\bar{u}$}
				& {$\bar{v}$}
				& {$\bar{p}$} \\
				
				\midrule
				
				{64}
				& \num{4.1820290702166904e-04}
				& \num{4.4787163267695810e-04}
				& \num{3.7895817316761381e-04}
				& \num{4.1820290702166904e-04}
				& {---}
				& {---}
				& {---}
				& {---} \\
				
				{128}
				& \num{5.2340139219698951e-05}
				& \num{5.6090569470117140e-05}
				& \num{4.7501955332684994e-05}
				& \num{5.2340139219698951e-05}
				& {3.00}
				& {3.00}
				& {3.00}
				& {3.00} \\
				
				{256}
				& \num{6.5403363307530014e-06}
				& \num{7.0105656671549904e-06}
				& \num{5.9393688898177827e-06}
				& \num{6.5403363307530014e-06}
				& {3.00}
				& {3.00}
				& {3.00}
				& {3.00} \\
				
				\bottomrule
				
			\end{tabular}
			
		\end{minipage}
	\end{center}
	
\end{table}

\begin{figure}[htb] 
	\centering
	\includegraphics[width=0.32\textwidth]{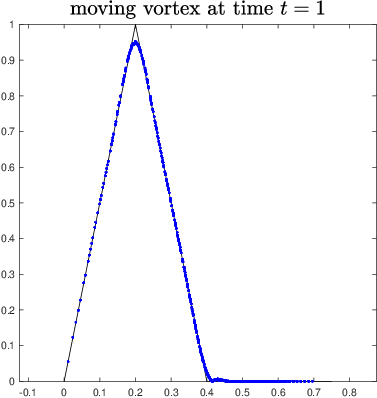}
	\hfill
	\includegraphics[width=0.32\textwidth]{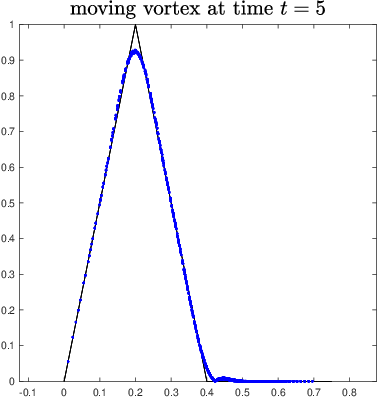}
	\hfill
	\includegraphics[width=0.32\textwidth]{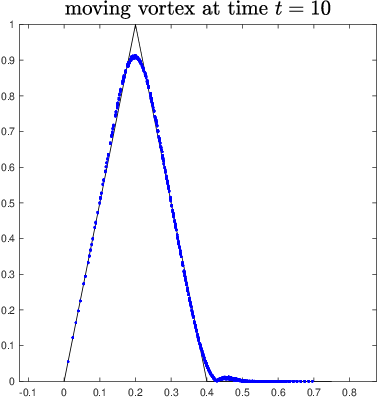}
	\includegraphics[width=0.32\textwidth]{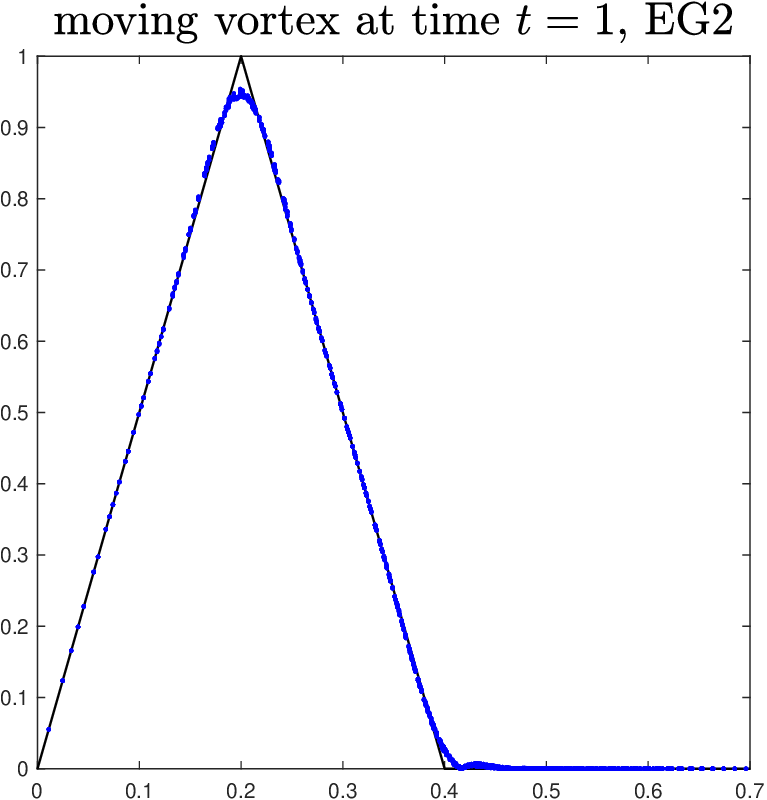}
	\hfill
	\includegraphics[width=0.32\textwidth]{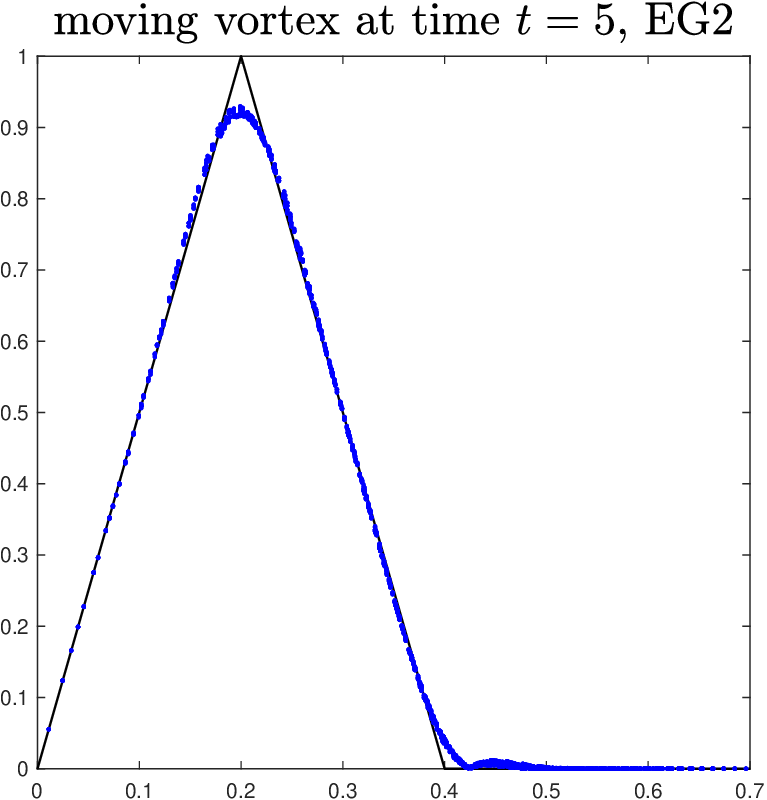}
	\hfill
	\includegraphics[width=0.32\textwidth]{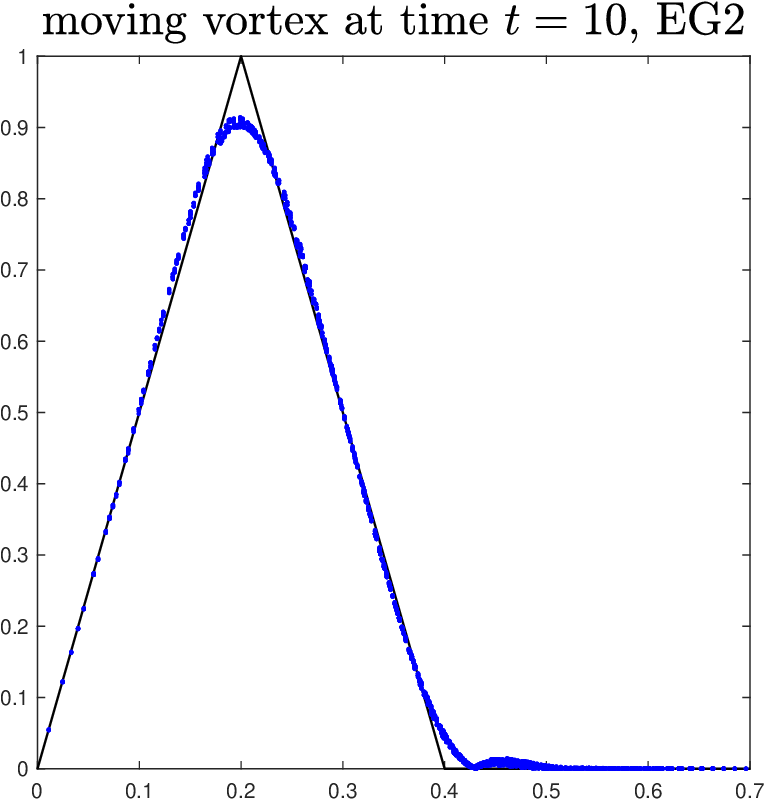}
	\caption{\label{fig:movingVortex}Moving vortex problem for the linearized Euler
		equations. Numerical results after one, five, and ten
		rotations using a grid with $64^2$ cells; Active Flux
		with moving coordinates approach for the update of point values
		(top) and with the EG2 operator (bottom). }
\end{figure}

Next we consider a moving vortex problem. The initial values 
in velocity and pressure agree with those of the vortex problem proposed for
acoustics by
Barsukow et al. \cite{article:BHKR2019}. By considering the same
data for the linearized Euler equations, we can study how well a moving
vortex is preserved.
\begin{example}
	We consider the linearized Euler equations with $u'=v'=1$,
	$\rho'=1$, $p'=1/1.4$, $c'=1$ on the domain $[-0.5,0.5]\times
	[-0.5,0.5]$ with periodic boundary conditions in the $x$- and $y$-
	directions. The initial conditions are
	\begin{equation}
		\begin{split}
			\rho(r,0) & = 0, \\
			\mathbf{u}(r,0)  & = 
			\mathbf{n} \left\{ \begin{array}{c c c}
				5r & : & 0\le r \le 0.2,\\
				2-5r & : & 0.2 < r \le 0.4 ,\\
				0 & : & r>0.4,
			\end{array}
			\right.\\
			p(r,0) & = 0.\\                    
		\end{split}
	\end{equation}
	At times $t=n$, $n\in \mathbb{N}$, the exact solution agrees with the
	initial condition. 
\end{example}
In Figure \ref{fig:movingVortex}, we show scatter plots of $\vert\mathbf{u}\vert$ at
times $t=1,5,10$, i.e., after one, five and ten rotations. We used a
relatively coarse mesh with only $64^2$ grid cells. The transported
vortex is no longer preserved by the Active Flux method with the exact
evolution operator. Both versions of the Active Flux method lead to
comparable and
accurate approximations.  

\begin{figure}[htb]
	\centerline{
		\includegraphics[width=0.32\textwidth]{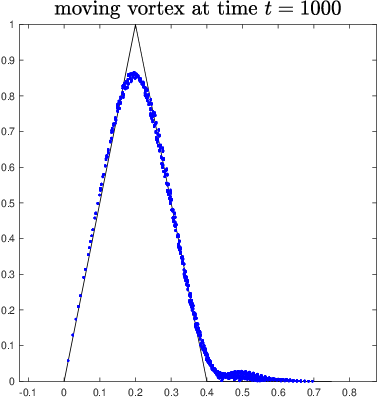}\hfil \hfil
		\includegraphics[width=0.32\textwidth]{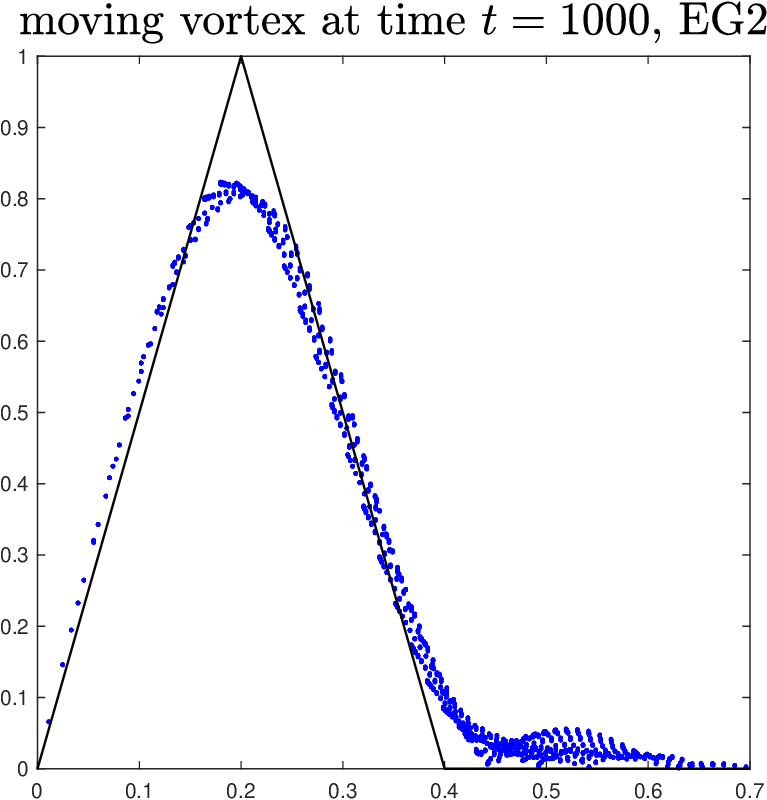}}
	\caption{\label{fig:movingVortex-2}Moving vortex problem for the linearized Euler
		equations. Numerical results after one rotations with
		$\vert\mathbf{u}'\vert/c' = 10^{-3}$ on a grid with $64^2$ cells.
		Results using the moving-grid approach (left) and the Active Flux
		method with the EG2 evolution operator (right). }
\end{figure}

More differences are visible when considering long-time approximations
in the regime $\sqrt{u'^2+v'^2}/c' \ll 1$. In Figure
\ref{fig:movingVortex-2}, we show numerical results of $\vert\mathbf{u}\vert$
at time $t=1000$ using $c'=1$ and $u'=v'=10^{-3}$. We use again a
coarse mesh with $64^2$ cells. In this flow regime, the new moving
grid approach leads to more accurate results.

\section{Fully discrete Active Flux methods for the  Euler equations}
\label{sec:5}
In this Section, we present numerical results for the fully discrete
Active Flux method for the Euler equations using both the approximate
and the exact evolution operator for the evolution of the
locally linearized problem.
The point value limiting strategy described in \cite{article:CHP2025}, and briefly reviewed in Section \ref{sec:limiting}, was
applied as needed to approximate test problems with shock waves.

First, we perform a convergence study on a smooth moving vortex in
order to compare the accuracy of the two third-order methods. 
\begin{example}\label{ex:movingVortexEuler}
	We consider a smooth rotating and traveling vortex, which propagates with constant velocity $(u_c,v_c)$. The initial values are given by
	\begin{align*}
		\rho(x,y,t) &=
		\begin{cases}
			\rho_c + \dfrac{1}{2}\left(1 - r^2\right)^6, & r < 1, \\
			\rho_c, & \text{otherwise},
		\end{cases} \\[1ex]
		u(x,y,t) &=
		\begin{cases}
			u_c - 1024 \sin(\theta)\,(1 - r)^6 r^6, & r < 1, \\
			u_c, & \text{otherwise},
		\end{cases} \\[1ex]
		v(x,y,t) &=
		\begin{cases}
			v_c + 1024 \cos(\theta)\,(1 - r)^6 r^6, & r < 1, \\
			v_c, & \text{otherwise},
		\end{cases} \\[1ex]
		p(x,y,t) &=
		\begin{cases}
			p_c + \bigl(p(r) - p(1)\bigr), & r < 1, \\
			p_c, & \text{otherwise}.
		\end{cases}
	\end{align*}
	Here, $r$ is the scaled distance from the initial center of the vortex, i.e.,
	\begin{align*}
		r = \sqrt{(x-0.5)^2 + (y-0.5)^2}/R.
	\end{align*}
	The function $p(r)$ is a polynomial of degree 36
	defined in \cite{article:CHP2025}. In our computations
	we use $R = 0.4, \rho_c = 0.5, u_c = v_c = 1$ and $p_c
	= 0.1. $ We simulate the vortex on the domain
	$\left[0, 1\right]^2$ using periodic boundary
	conditions in $x$- and $y$-directions. Originally, the
	test problem was proposed in \cite{article:KRM2008}.	
\end{example}

\begin{table}[!ht]
	\caption{\label{table:EuConMoving} Error measured in the $L_1$-norm and EOC for smooth moving
		vortex using moving-grid approach with $\mbox{CFL} = 0.45$ at $t=1$. }
	
	\sisetup{
		scientific-notation = true,
		round-mode = places
	}
	
	\begin{center}
		\begin{minipage}{0.9\textwidth}
			\begin{tabular}{
					*1{S[table-column-width=0.5cm,table-text-alignment=center]}
					*4{S[table-column-width=1.2cm,table-text-alignment=center]}
					*4{S[table-column-width=0.5cm,table-text-alignment=center,round-precision=4]}
				}
				\toprule
				
				{Res.} &
				\multicolumn{4}{c}{Error in} &
				\multicolumn{4}{c}{EOC in} \\
				
				\midrule
				
				&
				{$\bar{\rho}$} &
				{$\bar{\rho u}$} &
				{$\bar{\rho v}$} &
				{$\bar{E}$}
				&
				{$\bar{\rho}$} &
				{$\bar{\rho u}$} &
				{$\bar{\rho v}$} &
				{$\bar{E}$}
				\\
				
				\midrule
				
				{32}  &
				\num{4.73213326e-04} &
				\num{8.34324165e-04} &
				\num{8.19320038e-04} &
				\num{1.17760056e-03} &
				{---} &
				{---} &
				{---} &
				{---}
				\\
				
				{64} &
				\num{7.68183108e-05} &
				\num{1.31043475e-04} &
				\num{1.30260993e-04} &
				\num{1.85471450e-04} &
				{2.62} &
				{2.67} &
				{2.65} &
				{2.67}
				\\
				
				{128} &
				\num{1.03580297e-05} &
				\num{1.74540903e-05} &
				\num{1.73909776e-05} &
				\num{2.46313624e-05} &
				{2.89} &
				{2.91} &
				{2.90} &
				{2.91}
				\\
				
				{256} &
				\num{1.31448644e-06} &
				\num{2.21165717e-06} &
				\num{2.20602364e-06} &
				\num{3.12331107e-06} &
				{2.98} &
				{2.98} &
				{2.98} &
				{2.98}
				\\
				
				{512} &
				\num{1.64609892e-07} &
				\num{2.77074153e-07} &
				\num{2.76584565e-07} &
				\num{3.91598008e-07} &
				{3.00} &
				{3.00} &
				{3.00} &
				{3.00}
				\\
				
				\bottomrule
			\end{tabular}
		\end{minipage}
	\end{center}
\end{table}

\begin{table}[!ht]
	\caption{\label{table:EuConEG2} Error measured in the $L_1$-norm and EOC for smooth moving vortex
		using EG2 evolution with $\mbox{CFL} = 0.279$ at $t=1$. }
	
	
	\sisetup{
		scientific-notation = true,
		round-mode = places
	}
	
	\begin{center}
		\begin{minipage}{0.9\textwidth}
			
			\begin{tabular}{
					*1{S[table-column-width=0.5cm,table-text-alignment=center]}
					*4{S[table-column-width=1.2cm,table-text-alignment=center]}
					*4{S[table-column-width=0.5cm,table-text-alignment=center, round-precision=4]}
				}
				
				\toprule
				
				{Res.} &
				\multicolumn{4}{c}{Error in} &
				\multicolumn{4}{c}{EOC in} \\
				\midrule
				
				&
				{$\bar{\rho}$} &
				{$\bar{\rho u}$} &
				{$\bar{\rho v}$} &
				{$\bar{E}$}
				&
				{$\bar{\rho}$} &
				{$\bar{\rho u}$} &
				{$\bar{\rho v}$} &
				{$\bar{E}$}
				\\
				
				\midrule
				
				{32} &
				\num{5.78057154e-04} &
				\num{9.96758141e-04} &
				\num{9.95770954e-04} &
				\num{1.43177243e-03} &
				{---} &
				{---} &
				{---} &
				{---}
				\\
				
				{64} &
				\num{9.43933467e-05} &
				\num{1.60598767e-04} &
				\num{1.63708122e-04} &
				\num{2.34249894e-04} &
				{2.61} &
				{2.63} &
				{2.60} &
				{2.61}
				\\
				
				{128} &
				\num{1.28056525e-05} &
				\num{2.15412502e-05} &
				\num{2.21045884e-05} &
				\num{3.13976525e-05} &
				{2.88} &
				{2.90} &
				{2.89} &
				{2.90}
				\\
				
				{256} &
				\num{1.62647545e-06} &
				\num{2.73376687e-06} &
				\num{2.81378995e-06} &
				\num{3.99066618e-06} &
				{2.98} &
				{2.98} &
				{2.97} &
				{2.98}
				\\
				
				{512} &
				\num{2.03890945e-07} &
				\num{3.42907710e-07} &
				\num{3.53619307e-07} &
				\num{5.01228872e-07} &
				{3.00} &
				{3.00} &
				{2.99} &
				{2.99}
				\\
				
				\bottomrule
				
			\end{tabular}
		\end{minipage}
	\end{center}
\end{table}
In Tables \ref{table:EuConMoving} and \ref{table:EuConEG2}, we show the
results of a numerical convergence study for Example \ref{ex:movingVortexEuler}
using both the moving-grid approach as well as the EG2 evolution operator.
The numerical approximation of the cell average values on the next
finer grid is used as the reference solution for the error
computation.
Both methods converge with third-order. For the considered
CFL numbers, the measured error
is slightly smaller for the moving-grid
approach with the exact acoustic solver. 

Next we consider different two-dimensional Riemann problems, as
described in \cite{article:Schultz-Rinne1993,article:SRCG1993,lax1998solution}.
\begin{example}
	We consider classical two-dimensional Riemann
	problems.
	In the computational domain $[0,1]^2$, the initial values in
	primitive variables $\mathbf{u}=(\rho,u,v,p)$ are given by
	\begin{equation*}
		\mathbf{u}(x,y,0)=
		\begin{cases}
			\mathbf{u}_1, & \quad x > x_0,\ y > y_0,\\
			\mathbf{u}_2, & \quad x \le x_0,\ y > y_0,\\
			\mathbf{u}_3, & \quad x \le x_0,\ y \le y_0,\\
			\mathbf{u}_4, & \quad x > x_0,\ y \le y_0,
		\end{cases}
	\end{equation*}
	where $\mathbf{u}_1,\dots,\mathbf{u}_4$ denote
	constant states
	prescribed in the four quadrants.
	The initial condition is therefore composed of four
	constant states
	separated by discontinuities aligned with the
	coordinate directions,
	intersecting at $(x_0,y_0)$. We use outflow boundary
	conditions.
	
	We will consider the following four configurations from
	\cite{lax1998solution}:

\begin{alignat*}{2}
	&\text{Configuration 3:}\qquad x_0=y_0=0.8, \\
	&\mathbf{u}_1 = (1.5,0,0,1.5),      &\quad
	&\mathbf{u}_2 = (0.5323,1.206,0,0.3),\\
	&\mathbf{u}_3 = (0.138,1.206,1.206,0.029), &\quad
	&\mathbf{u}_4 = (0.5323,0,1.206,0.3).\\[0.5em]
	&\text{Configuration 4:}\qquad x_0=y_0=0.5, \\
	&\mathbf{u}_1 = (1,1,0,0,1.1),   &\quad
	&\mathbf{u}_2 = (0.5065,0.8939,0,0.35),\\
	&\mathbf{u}_3 = (1.1,0.8939,0.8939,1.1),        &\quad
	&\mathbf{u}_4 = (0.5065,0,0.8939,0.35).\\[0.5em]
	&\text{Configuration 12:}\qquad x_0=y_0=0.5, \\
	&\mathbf{u}_1 = (0.5313,0,0,0.4),   &\quad
	&\mathbf{u}_2 = (1,0.7276,0,1),\\
	&\mathbf{u}_3 = (0.8,0,0,1),        &\quad
	&\mathbf{u}_4 = (1,0,0.7276,1).\\[0.5em]
	&\text{Configuration 17:}\qquad x_0=y_0=0.5, \\
	&\mathbf{u}_1 = (1,0,-0.4,1),       &\quad
	&\mathbf{u}_2 = (2,0,-0.3,1),\\
	&\mathbf{u}_3 = (1.0625,0,0.2145,0.4), &\quad
	&\mathbf{u}_4 = (0.5197,0,-1.1259,0.4).
\end{alignat*}

\end{example}

The first Riemann problem we consider is the relatively simple
configuration 12.
The initial data consist of two forward-moving
shocks located between the first and second quadrants, and between the first
and fourth quadrants. 
Stationary contact discontinuities separate the remaining quadrants.
The shock interaction produces a Mach reflection pattern, while the
interaction of the contact lines creates a vortex structure.  While this problem can be approximated by our fully discrete Active Flux method with the EG2 evolution operator without any limiting as shown in \cite[Figure 6]{article:CHP2025}, some inaccuracies near the thin Mach stem that forms between the two shocks  appeared when the unlimited Active Flux method using the moving grid approach was used. Applying the point value limiting strategy described in Section \ref{sec:limiting} and in more detail in \cite[Section 3]{article:CHP2025} suppresses this behavior.

\begin{figure}	
	\includegraphics[width=0.32\linewidth]{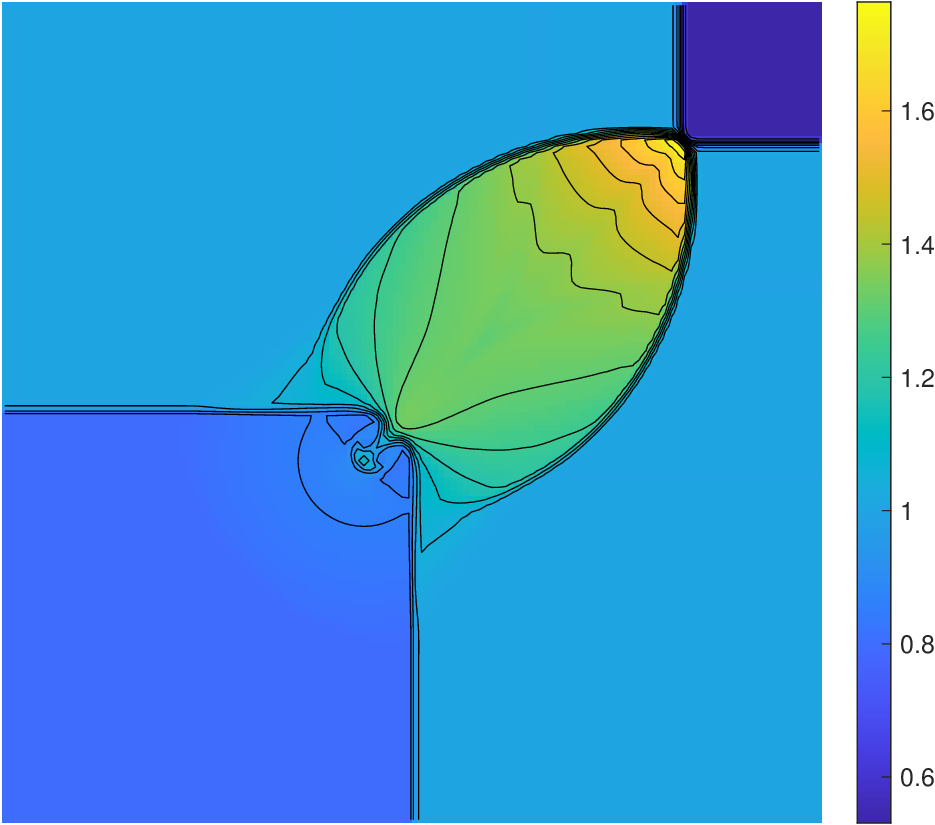}
	\includegraphics[width=0.32\linewidth]{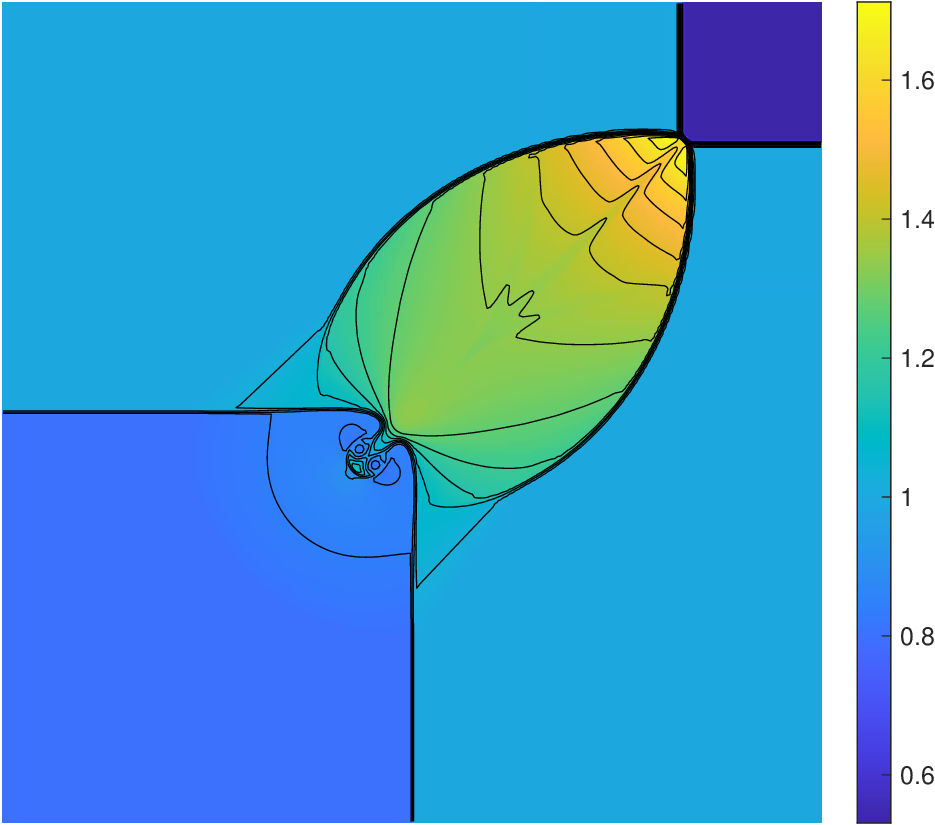}
	\includegraphics[width=0.32\linewidth]{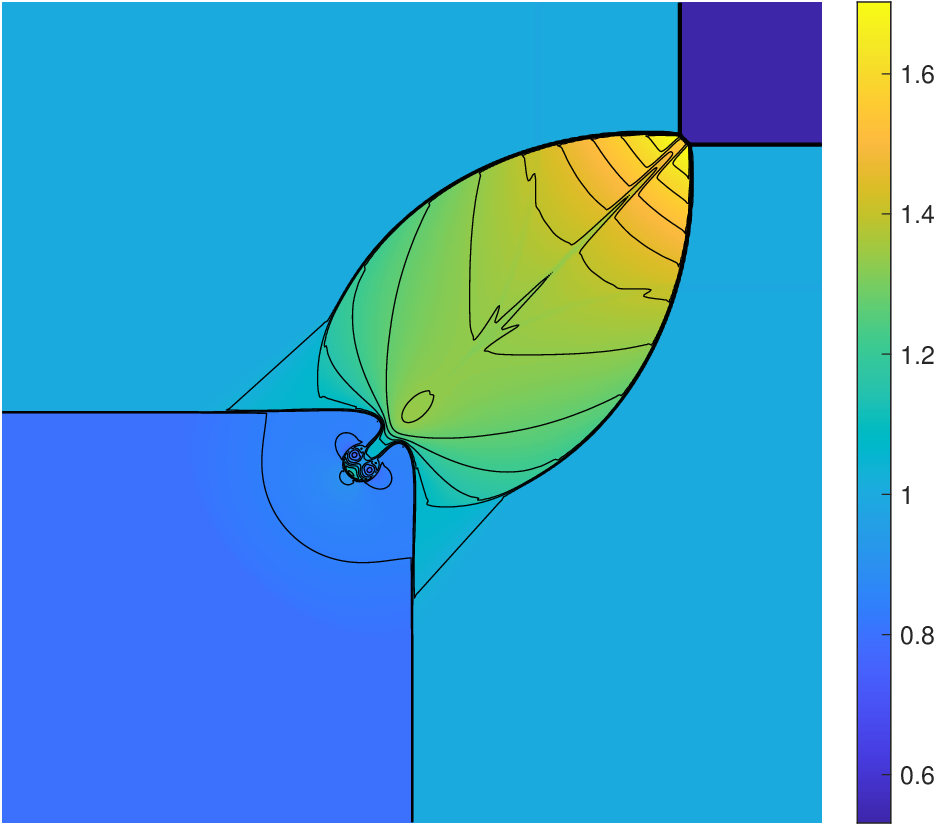} \\[0.15cm]
	\includegraphics[width=0.32\linewidth]{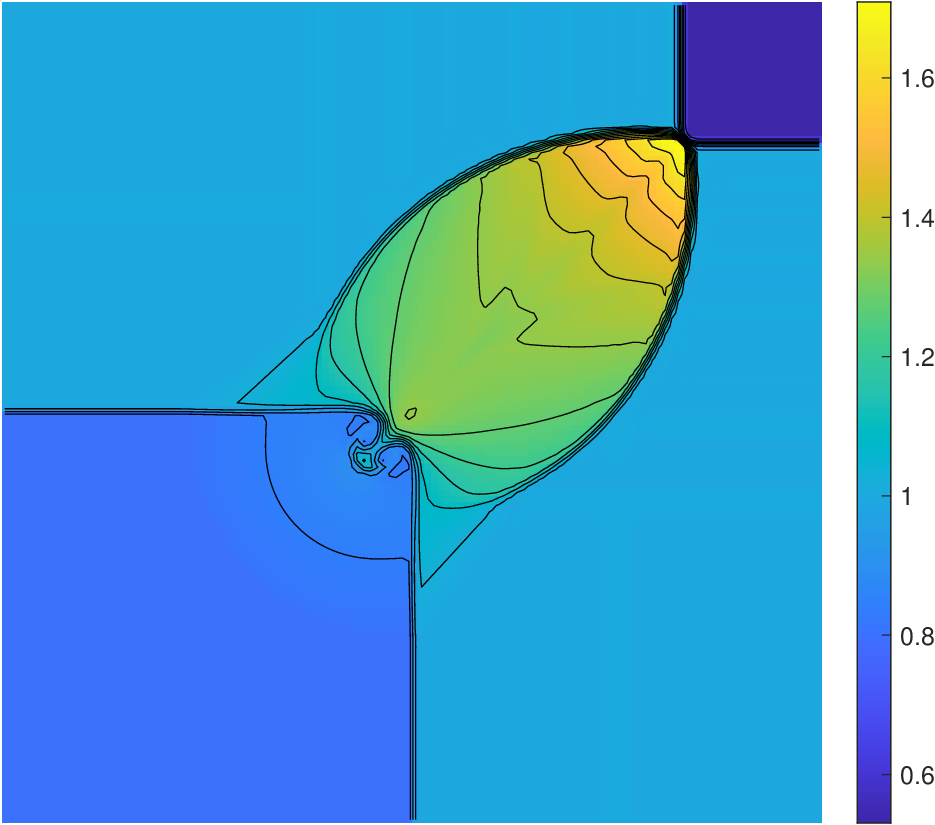}
	\includegraphics[width=0.32\linewidth]{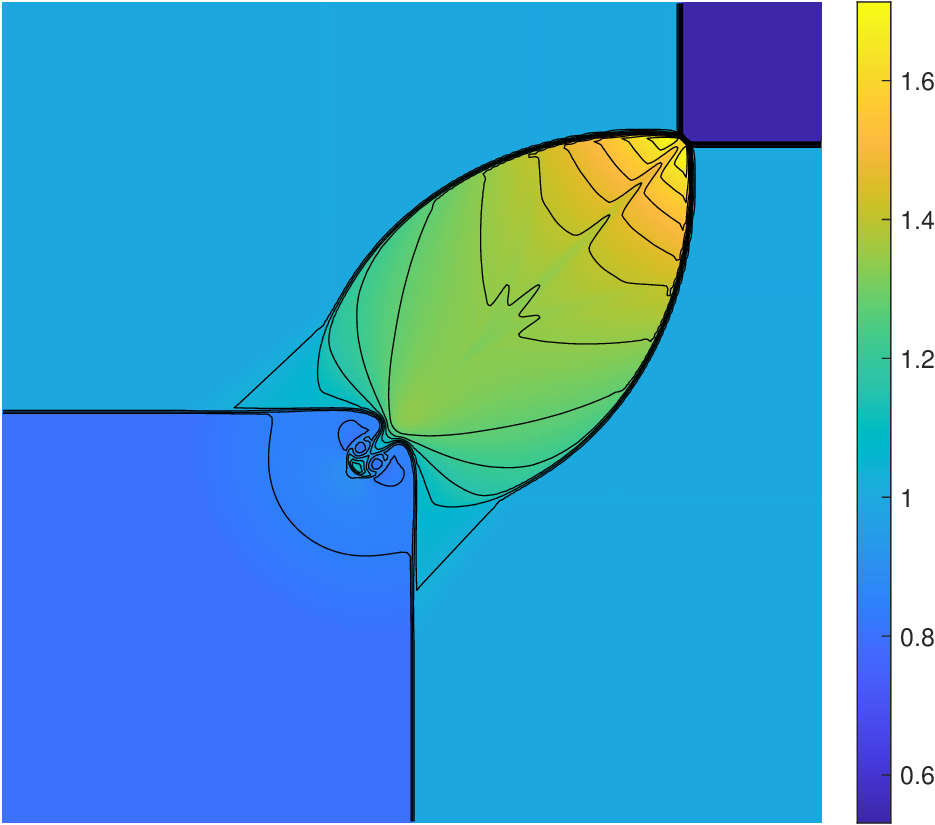}
	\includegraphics[width=0.32\linewidth]{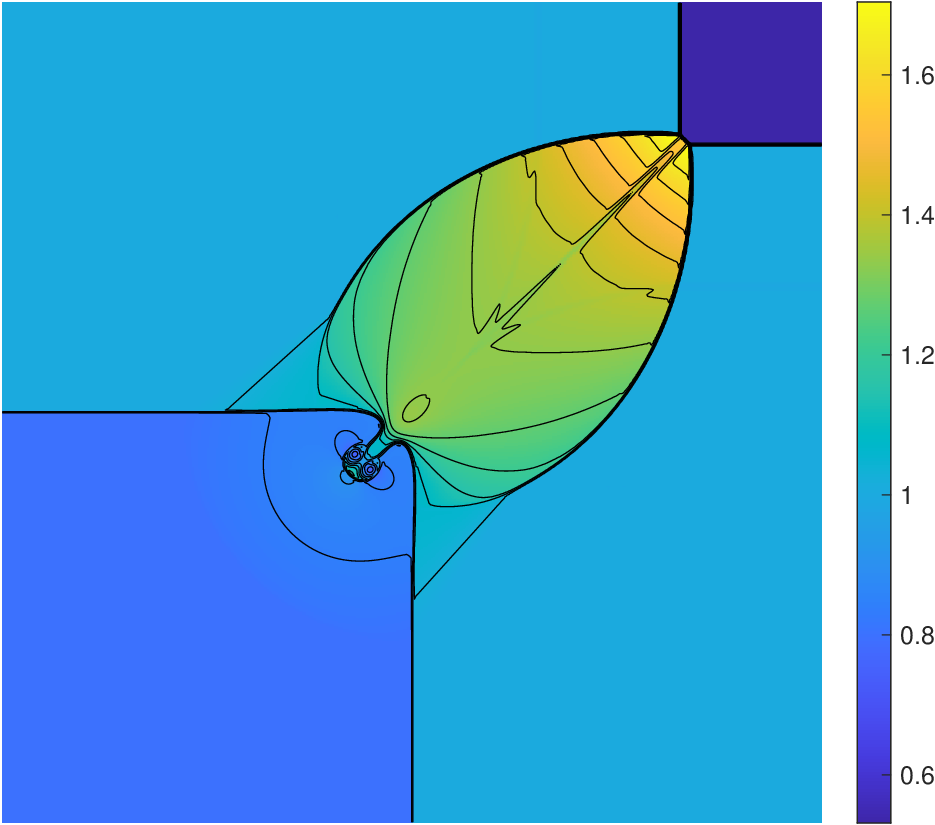}
	\caption{\label{fig:C12} Numerical results for
		Configuration 12 on grids with $128^2$ (left), $256^2$ (middle), and
		$512^2$ (right) grid cells. Plot of density at time $t=0.21$ using
		moving-grid approach (top) and EG2 (bottom).}
\end{figure}

Figure \ref{fig:C12} shows solutions computed using the exact evolution
operator and the approximate EG2 operator on different grids at time
$t=0.21$. For all computations, the limiting of the linearization and the point values was applied.
\ignore{For this test problem, the moving-grid approach required a limiting of the local linearization, as outlined in Section \ref{sec:limiting}. The unlimited version led to some inaccuracies near the thin Mach stem that forms between the two shocks. When using the EG2 evolution operator, we did not observe such problems. For the simulations shown in Figure \ref{fig:C12}, we used the same limited local linearization for both.}
Both versions of the method provide a good resolution of the structures even on
coarse grids, and the resulting solutions are very similar.

The next example, configuration 4,  consists of four moving shocks. As the solution evolves, curved shock fronts emerge.  For the method based on the moving grid approach,  the modified linearization described in Section \ref{sec:limiting} was needed. Without this modification, unphysical artefacts develop in the vicinity of the curved shock fronts. 
This observation suggests that configuration 4 may provide an interesting test case for other Active Flux methods based on exact evolution operators, such as those proposed in \cite{article:Barsukow2025} and \cite{duraisamy2026}.
Without limiting, both methods produce visible unphysical oscillations. Therefore, the limiting strategy for point values described in Section \ref{sec:limiting} was used for the computations. No additional flux limiting was needed.
\begin{figure}	
	\includegraphics[width=0.32\linewidth]{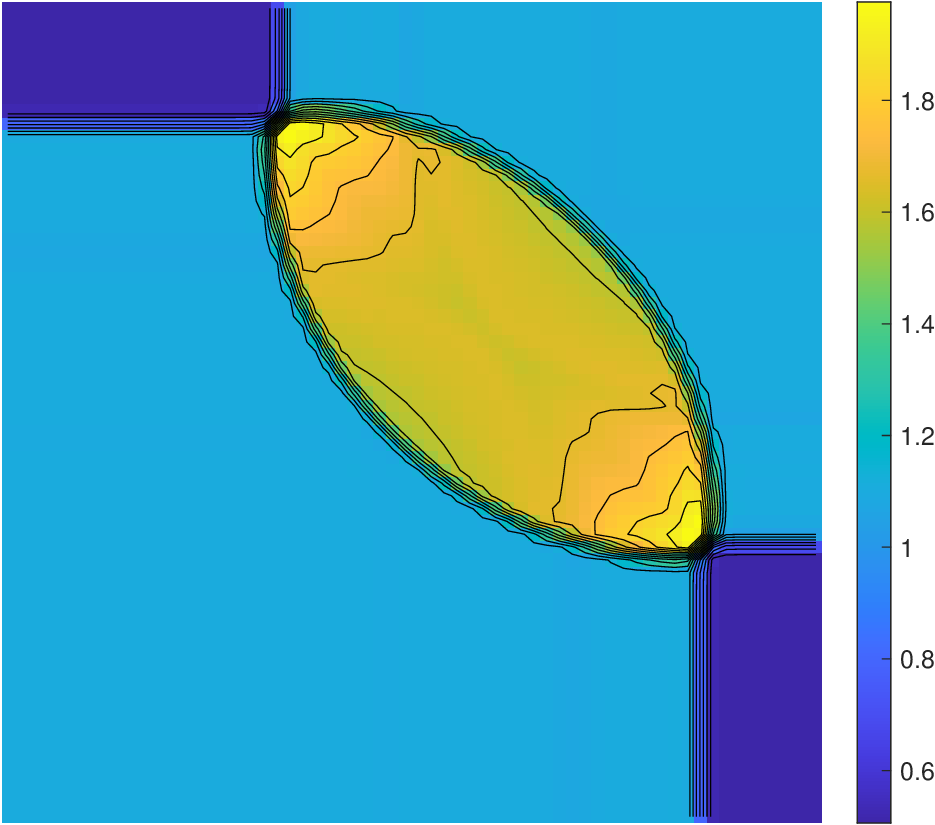}
	\includegraphics[width=0.32\linewidth]{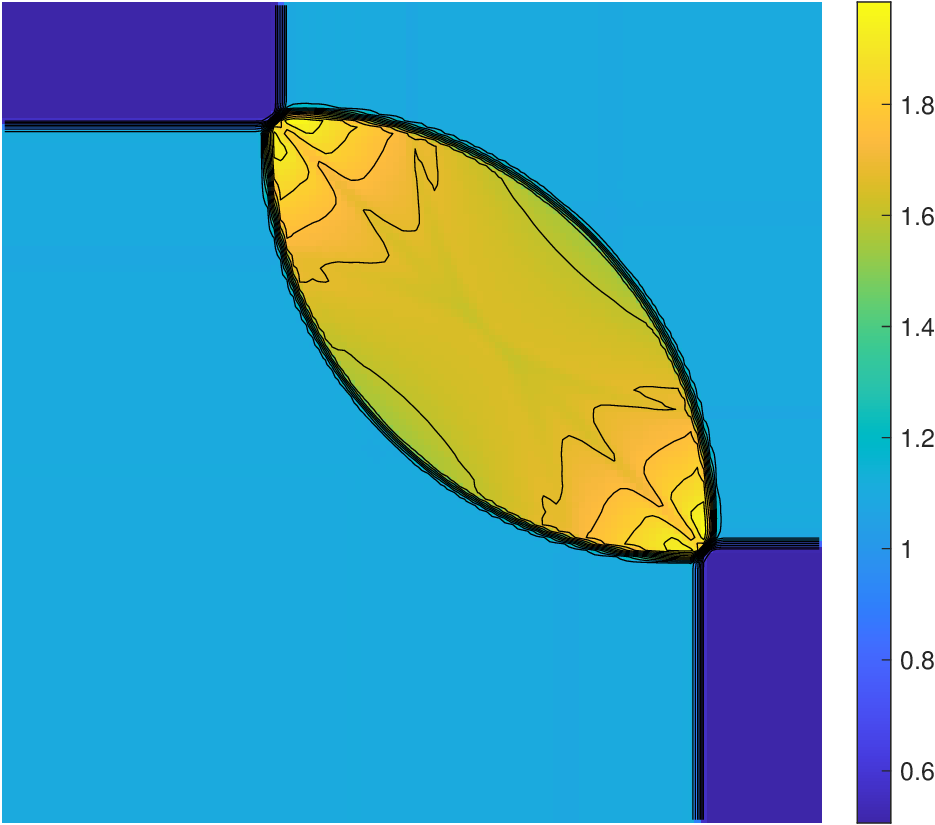}
	\includegraphics[width=0.32\linewidth]{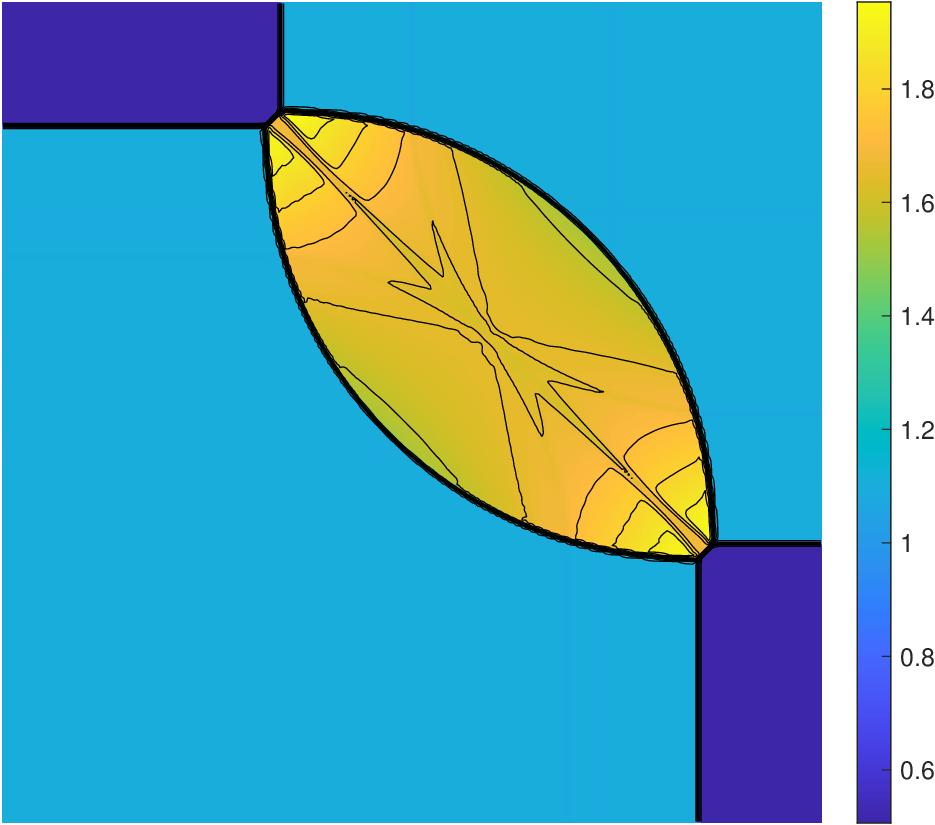} \\[0.15cm]
	\includegraphics[width=0.32\linewidth]{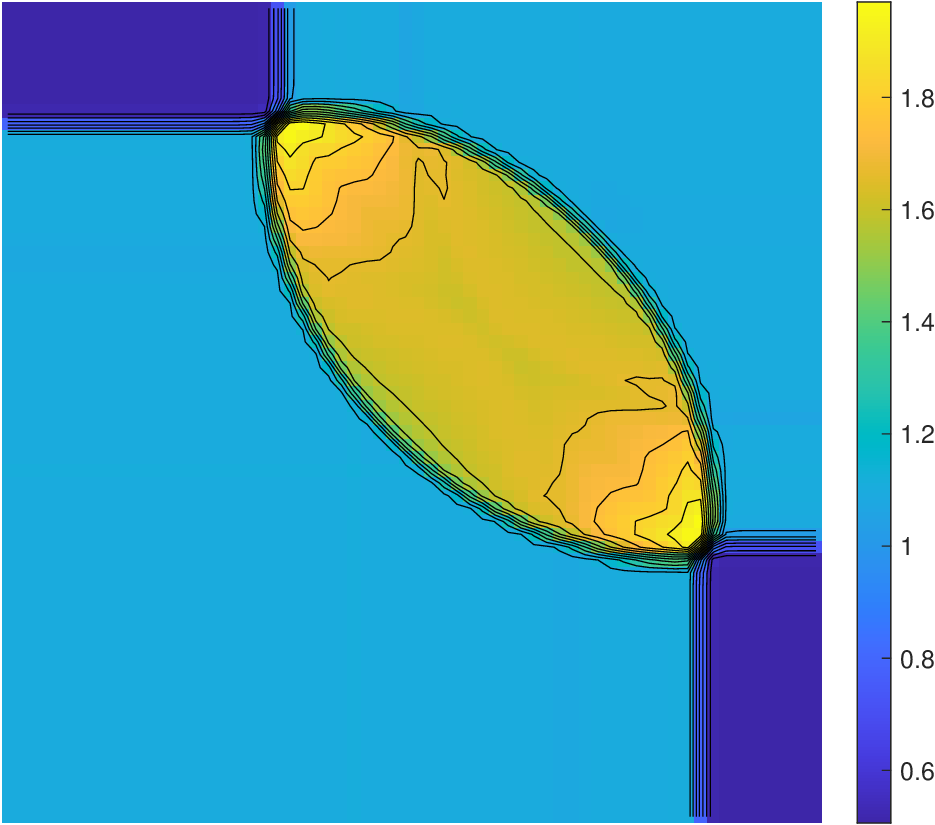}
	\includegraphics[width=0.32\linewidth]{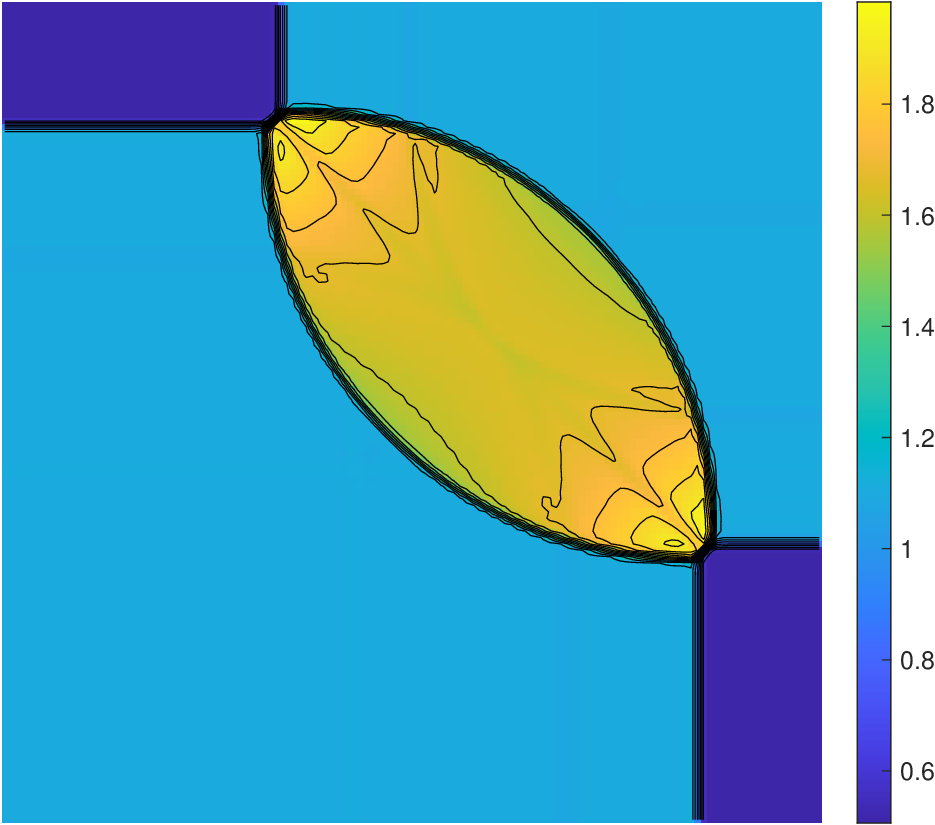}
	\includegraphics[width=0.32\linewidth]{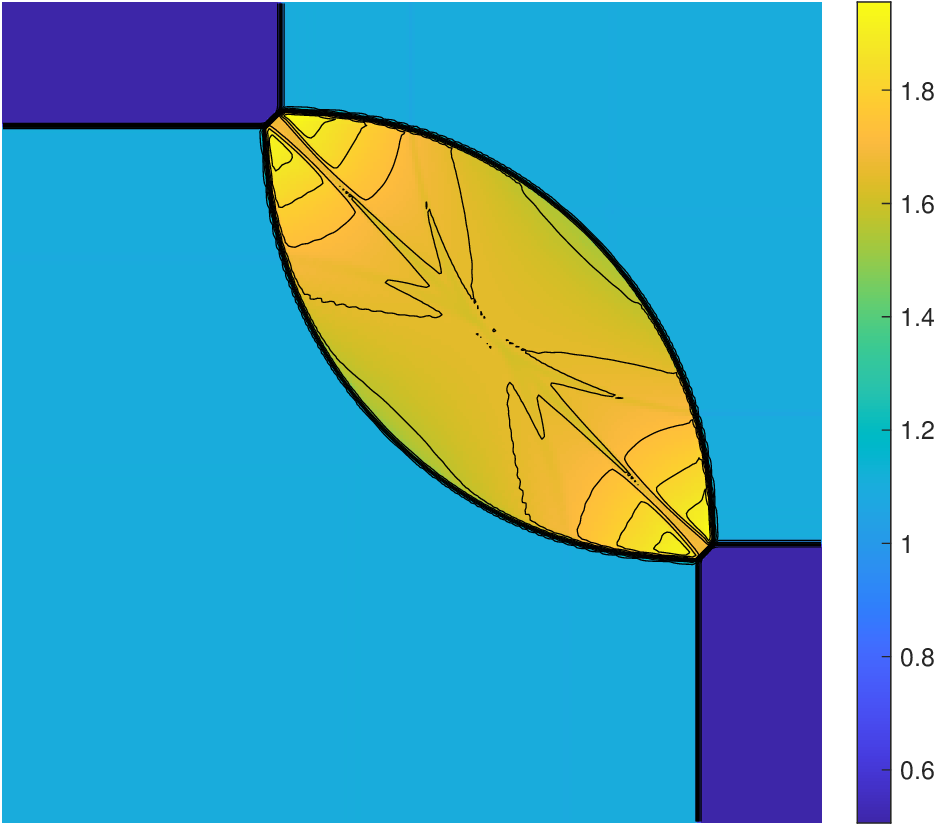}
	\caption{\label{fig:C4} Numerical results for
		Configuration 4 on grids with $64^2$ (left), $128^2$ (middle), and
		$256^2$ (right) grid cells. Plot of density at time $t=0.21$ using
		moving-grid approach (top) and EG2 (bottom).}
\end{figure}
Figure \ref{fig:C4} shows solutions computed with both versions of the method on different grids at time $t=0.21$. As in the previous example, only minor differences between the two methods can be observed. Note that the coarsest simulations were performed on grids with only  $64^2$ cells. 

\begin{figure}	
	\includegraphics[width=0.32\linewidth]{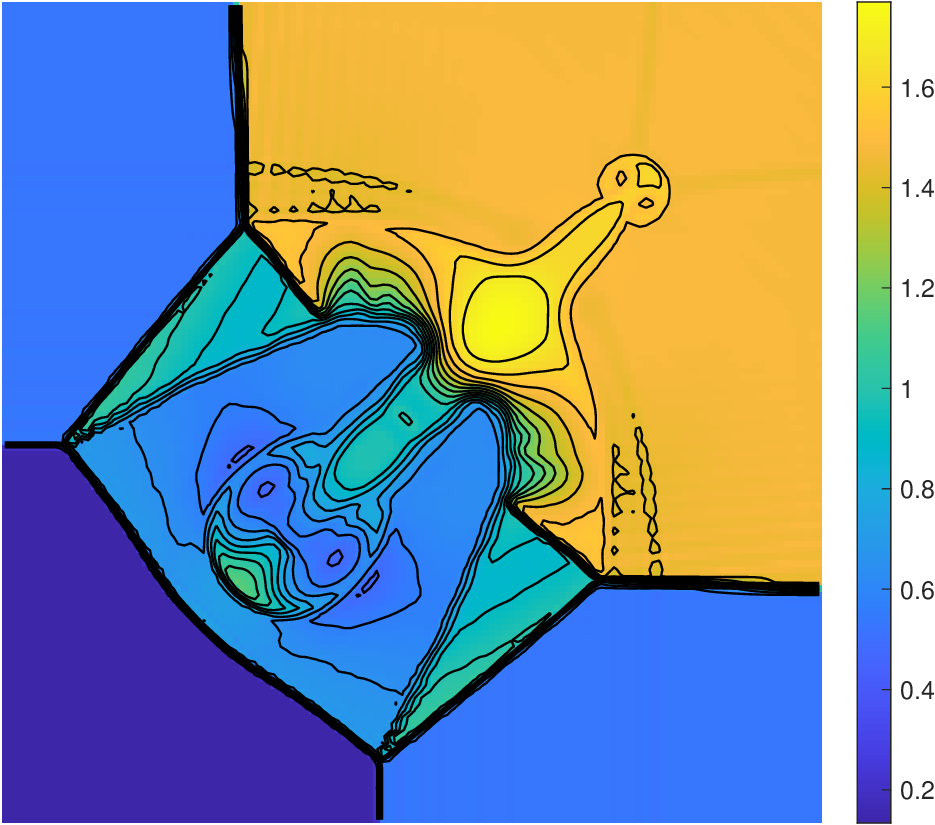}\hfill
	\includegraphics[width=0.32\linewidth]{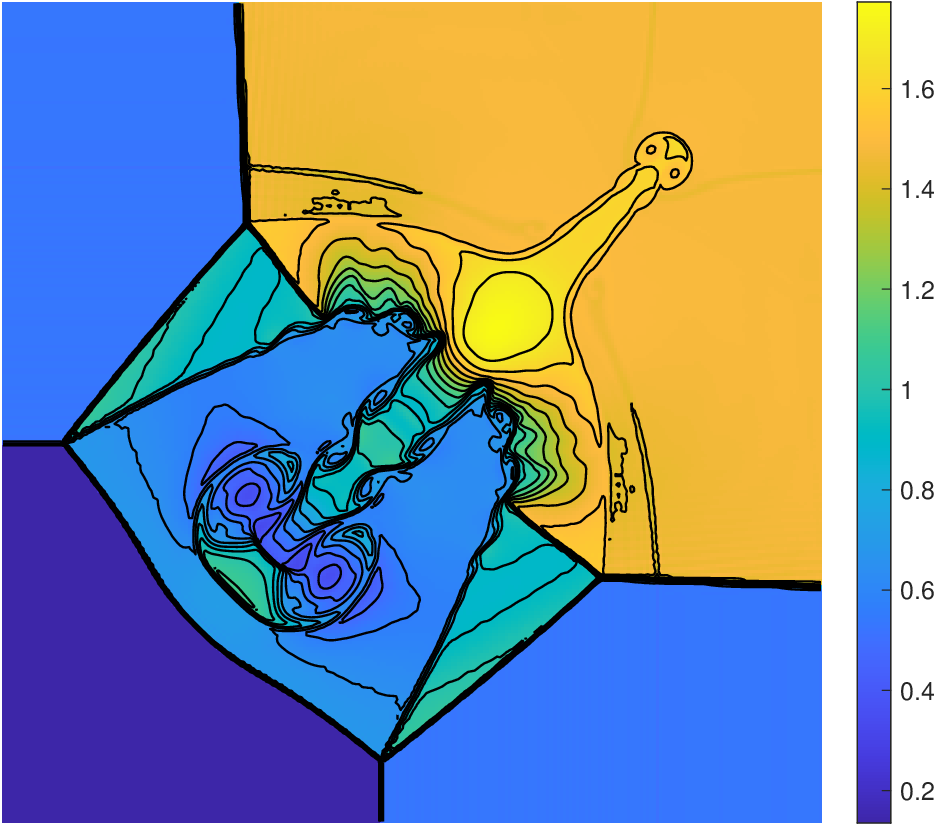}\hfill
	\includegraphics[width=0.32\linewidth]{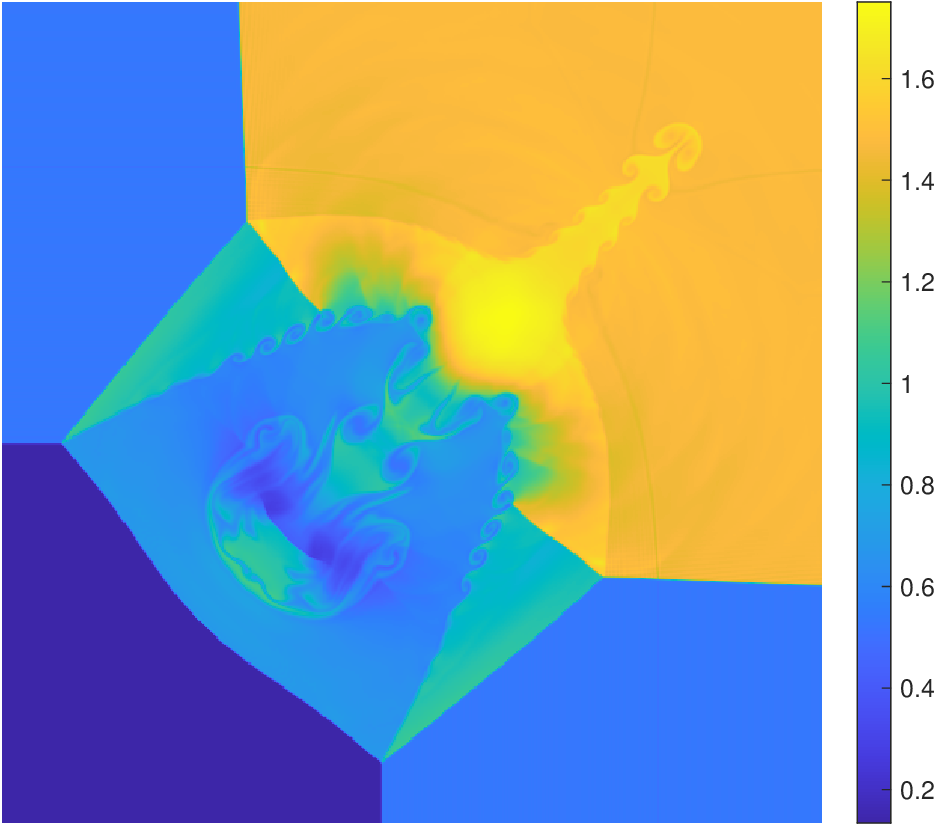} \\[0.15cm]
	\includegraphics[width=0.32\linewidth]{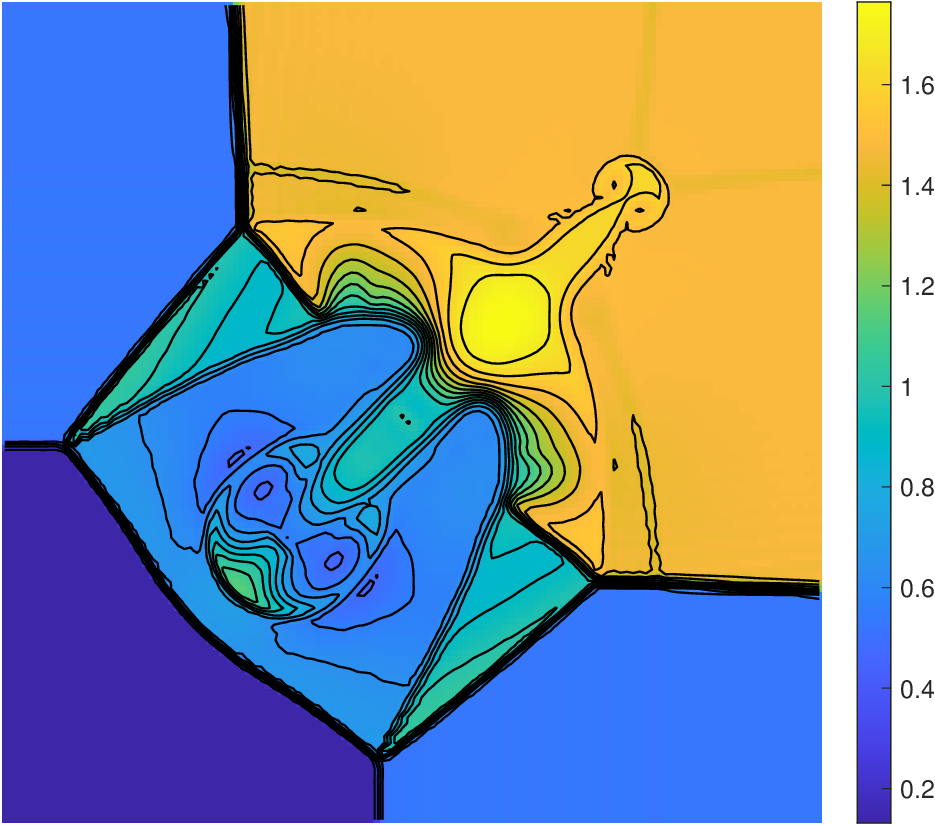}\hfill
	\includegraphics[width=0.32\linewidth]{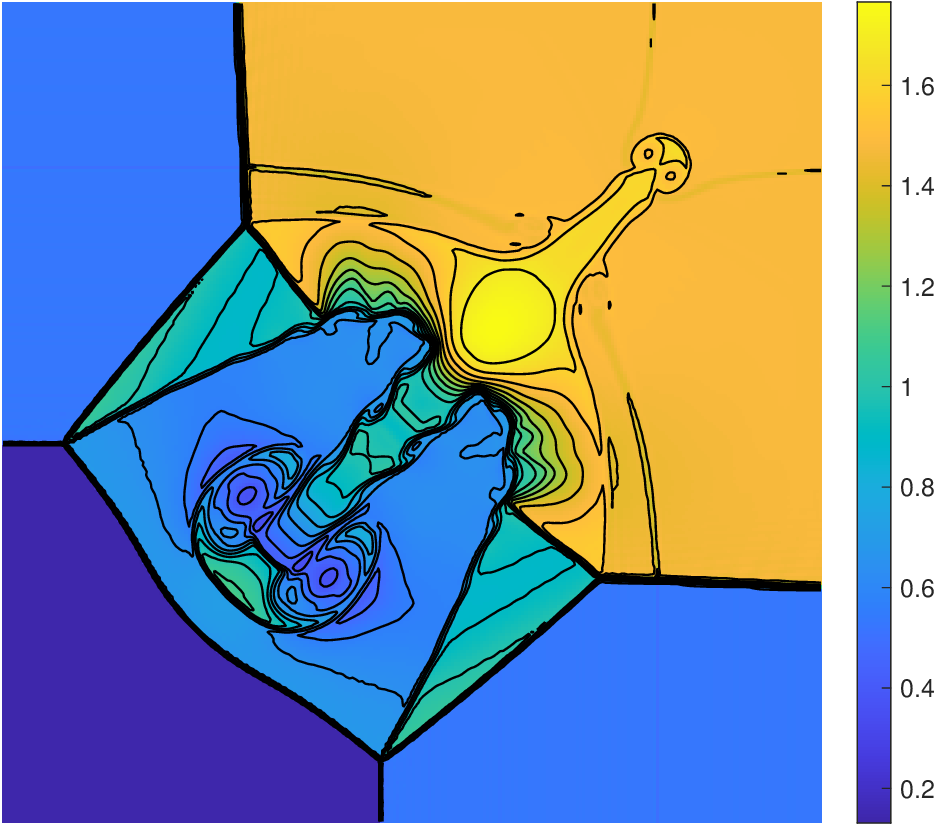}\hfill
	\includegraphics[width=0.32\linewidth]{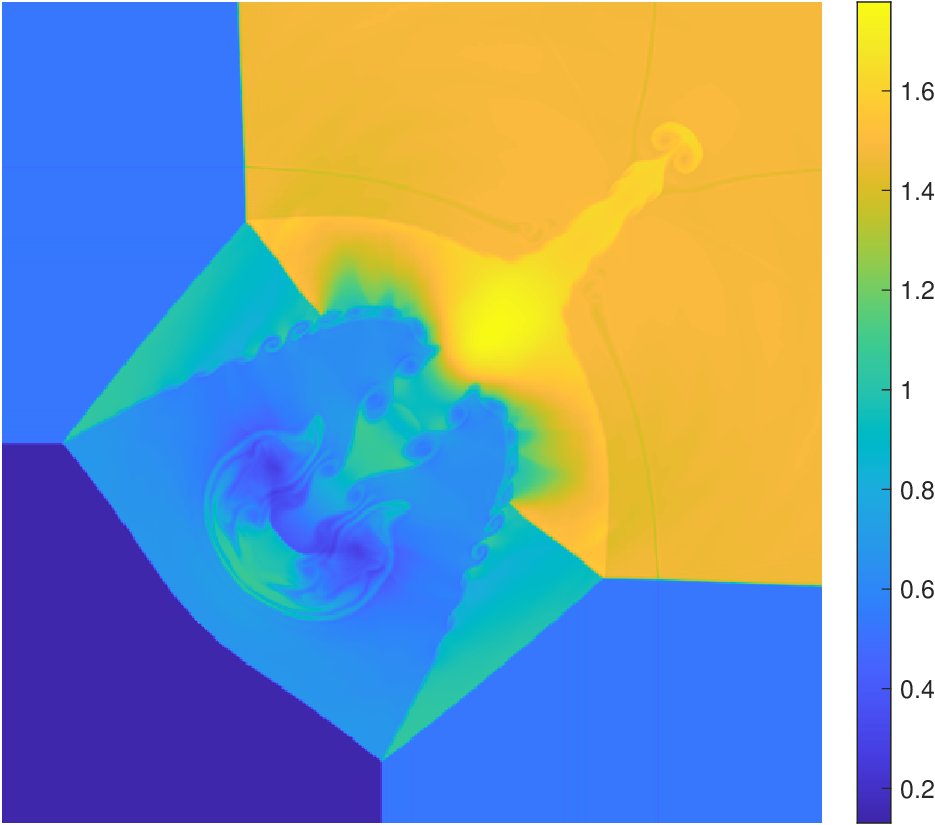}
	\caption{\label{fig:C3} Numerical results for
		Configuration 3 on grids with $128^2$ (left), $256^2$ (middle), and
		$512^2$ (right) grid cells. Plot of density at time $t=0.8$ using
		moving-grid approach (top) and EG2 (bottom).}
\end{figure}

The next example, i.e., configuration 3, 
consists of four moving shocks, whose interaction produces a
double Mach reflection pattern. 
For this example, both versions of the method require the limiting strategy for point values. 
No additional flux limiting to obtain bound
preserving cell average values was needed.
Over time, a complex solution structure develops, including contact discontinuities and vortical structures which eventually lead to a Kelvin-Helmholtz instability. 
Once again, both methods produce comparable approximations of the solution, although minor differences are visible.
Figure \ref{fig:C3} shows solutions obtained on grids with $128^2$, $256^2$, and $512^2$ cells at $t=0.8$.

In our next test problem, which corresponds to configuration
17, we consider a slight
modification of the standard definition of this two-dimensional
Riemann problem.
The solution structure consists of two shocks on the vertical
interface, interacting with two contact discontinuities on the
horizontal interface and producing a central shock-contact-line
interaction pattern with strong shear around the origin.
Highly resolved computations, compare with \cite{article:JN2018}, show a
Kelvin-Helmholtz instability along the stationary contact line between
the third and fourth quadrants. In our numerical simulations, using the moving-grid approach, we observe that the onset of the Kelvin-Helmholtz instability also
dependents on the definition of the initial values. If we set the initial point values along the lines $x=x_0$ and $y=y_0$ to be equal to the average of the two neighboring constant values, the numerical solution will be less susceptible to the Kelvin-Helmholtz instability. We did not observe this effect when using the EG2 evolution operator but used the averaged values along the lines $x=x_0$, $y=y_0$ to define the initial values for both simulations.

To trigger the instability even on coarse
grids, we further changed the initial values by adding a random perturbation
with amplitude $10^{-5}$ to the velocity components
of the point values along the line segment $x=x_0$, $y<y_0$.  
\begin{figure}	
	\includegraphics[width=0.32\linewidth]{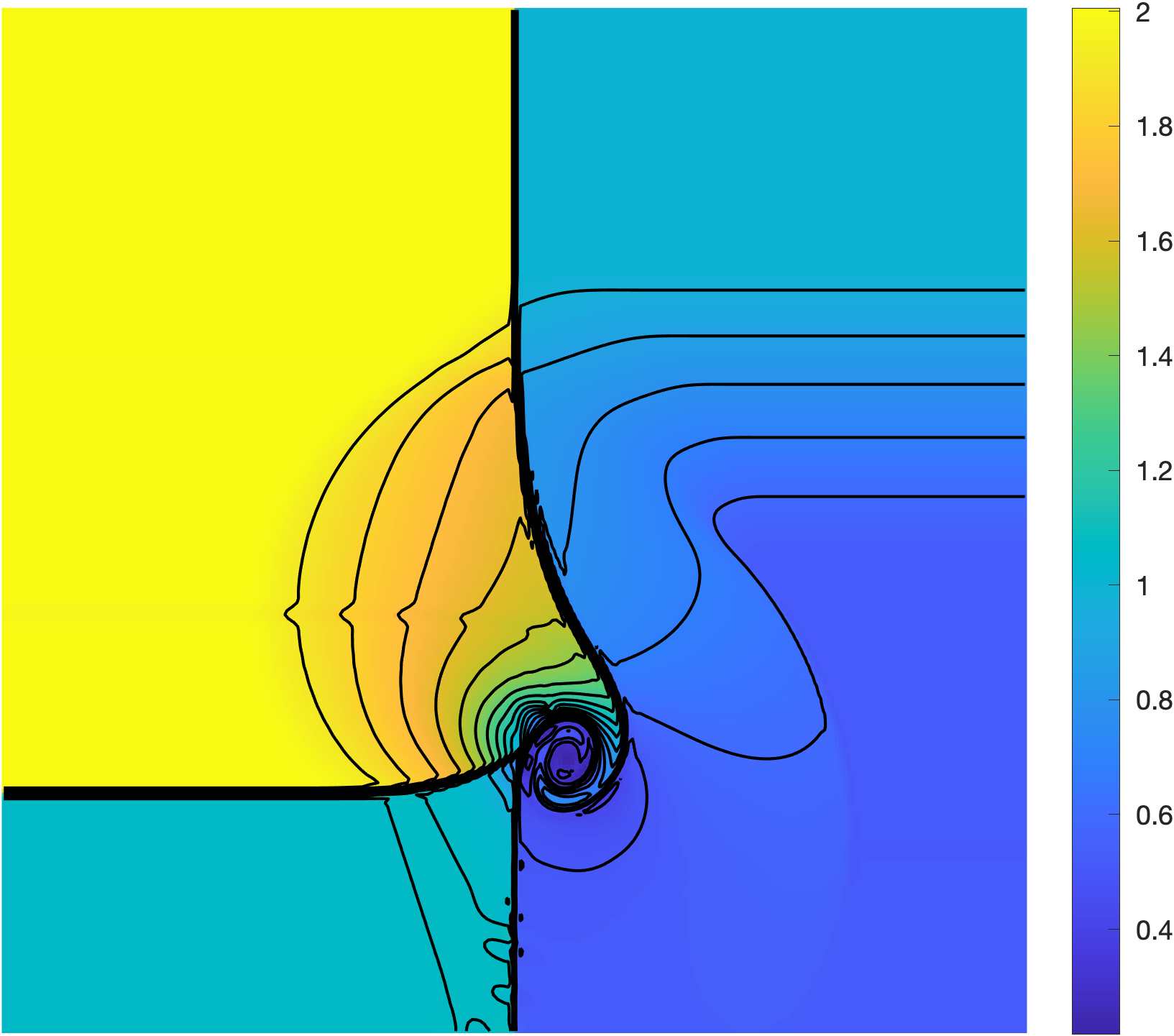}\hfill
	\includegraphics[width=0.32\linewidth]{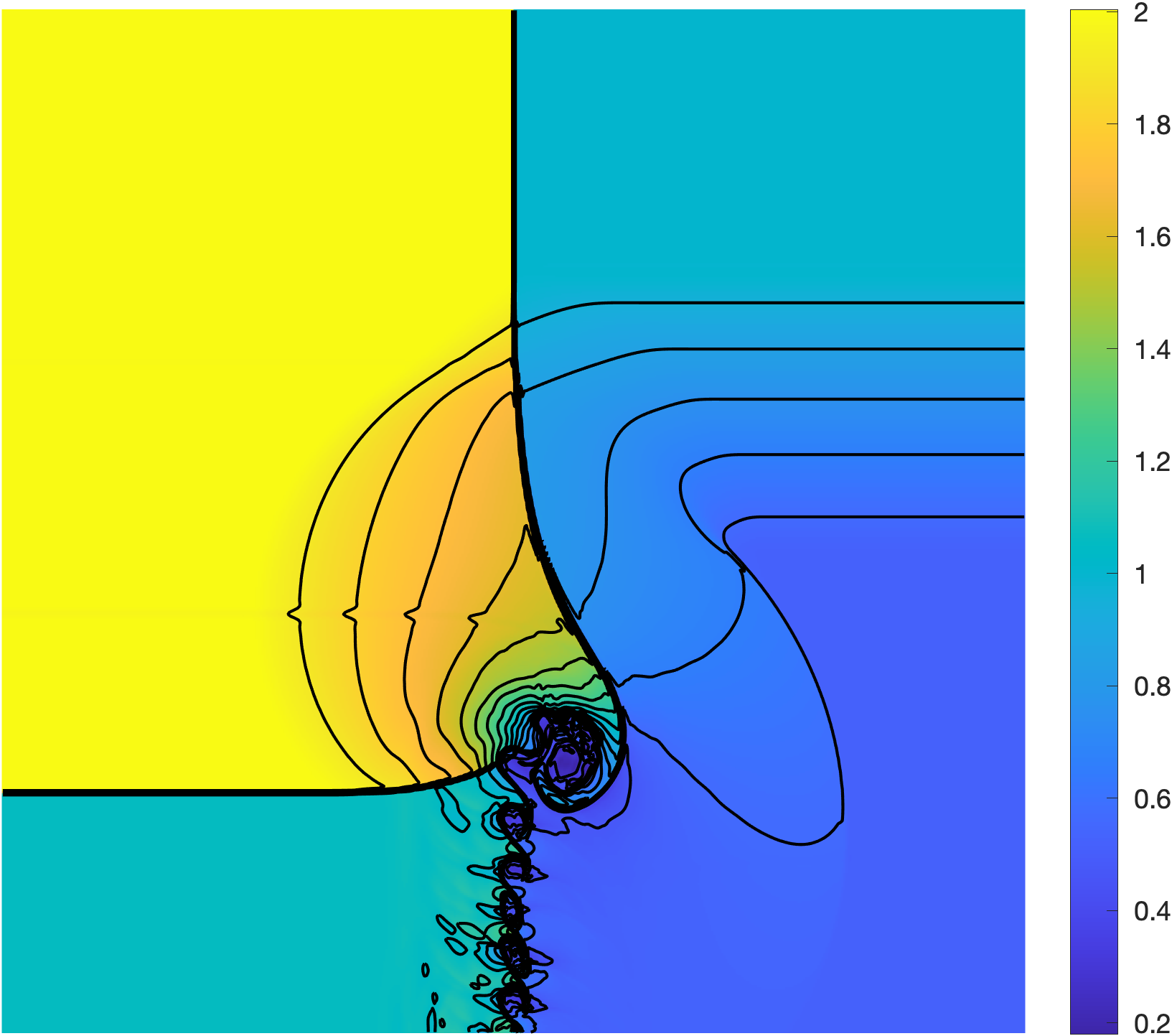}\hfill
	\includegraphics[width=0.32\linewidth]{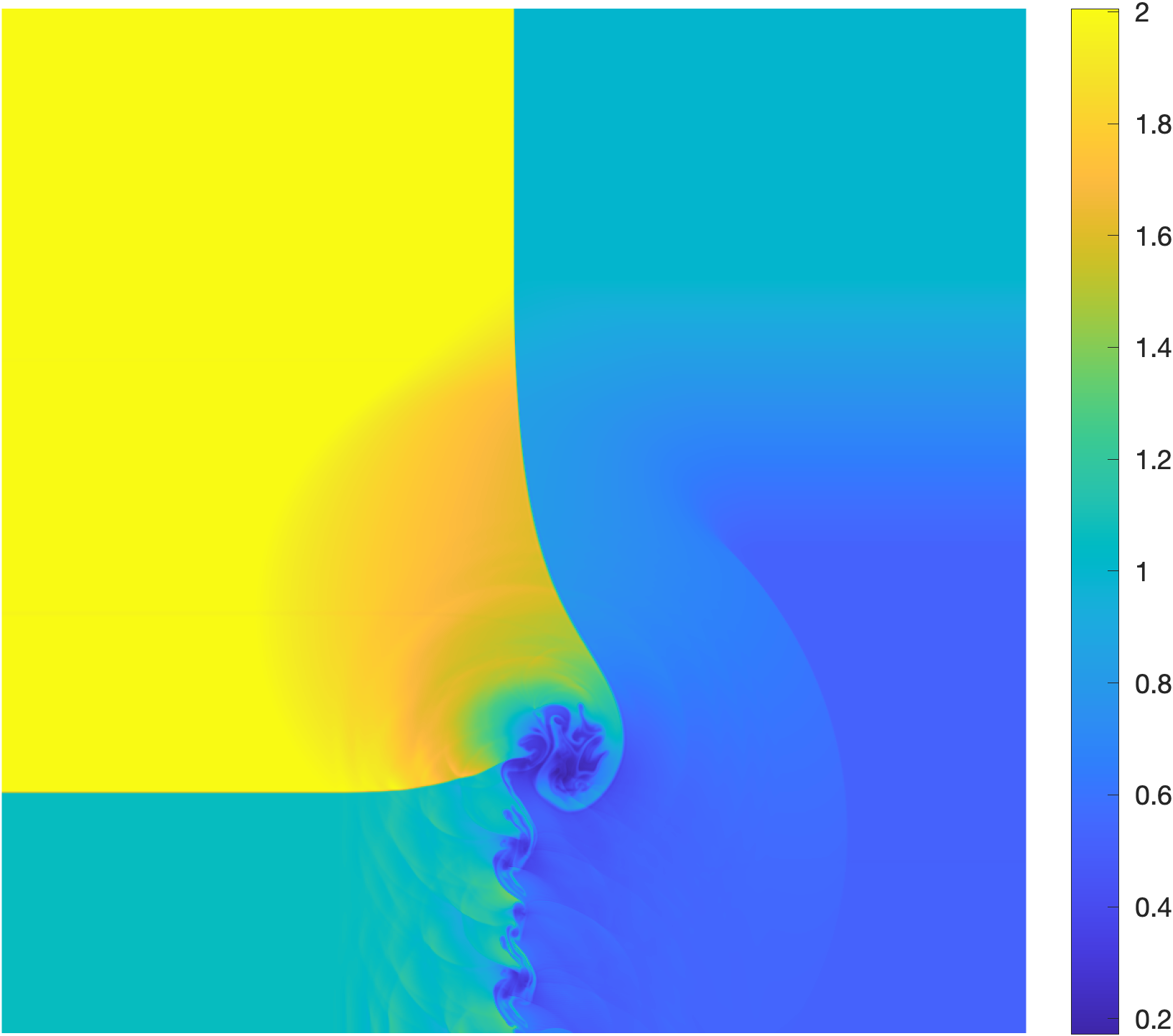}
	\\
	
	\includegraphics[width=0.32\linewidth]{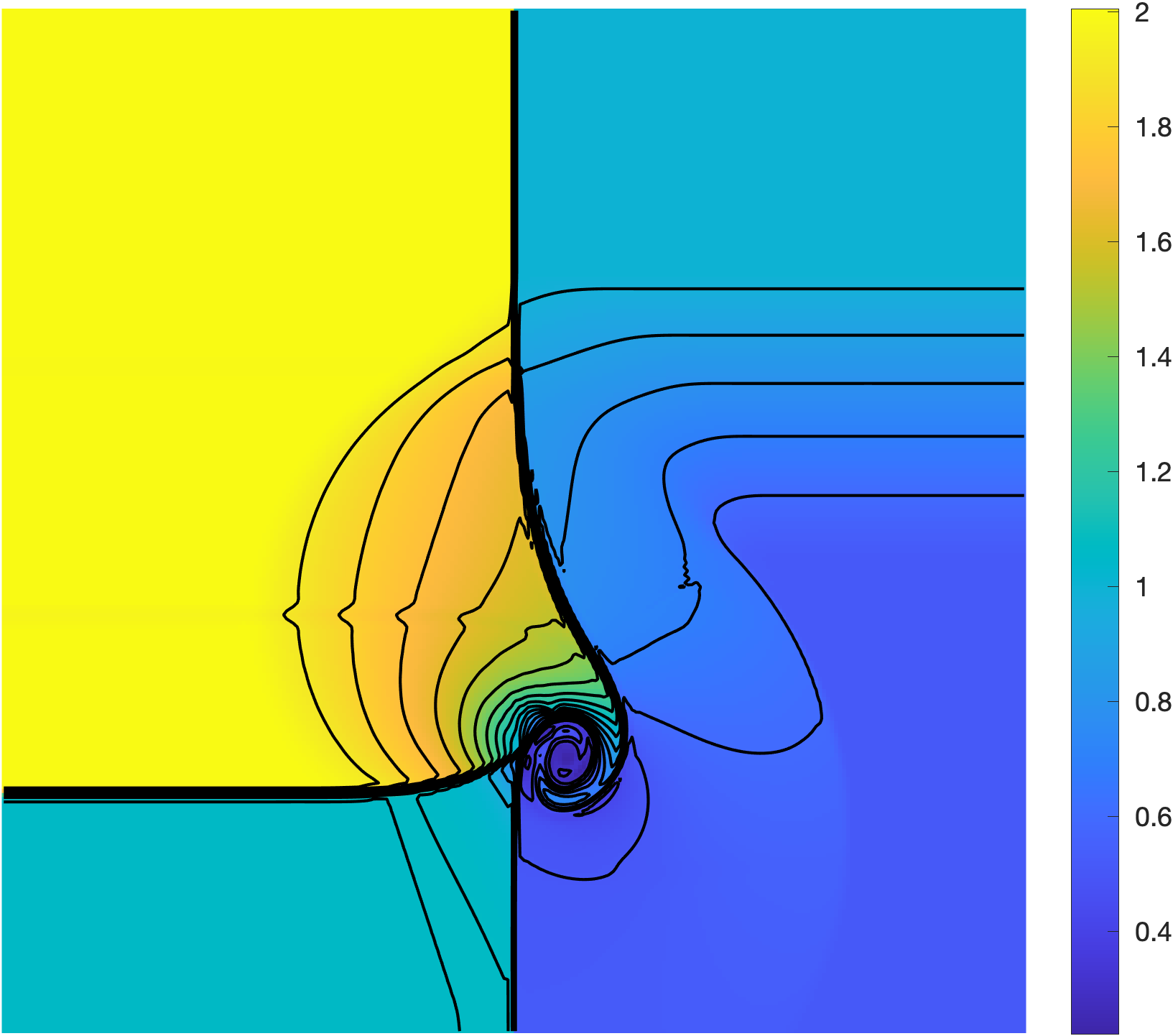}\hfill
	\includegraphics[width=0.32\linewidth]{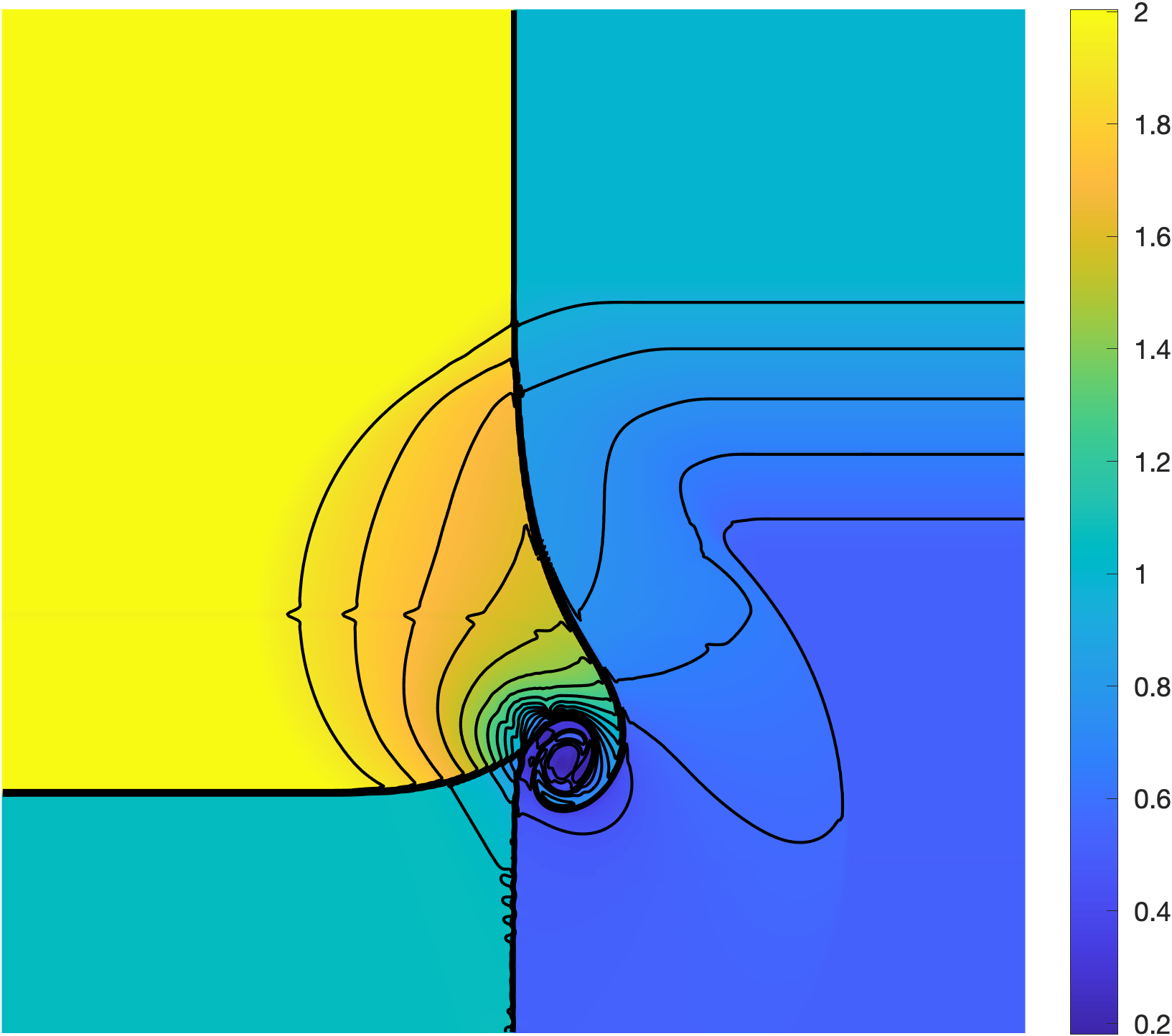}\hfill
	\includegraphics[width=0.32\linewidth]{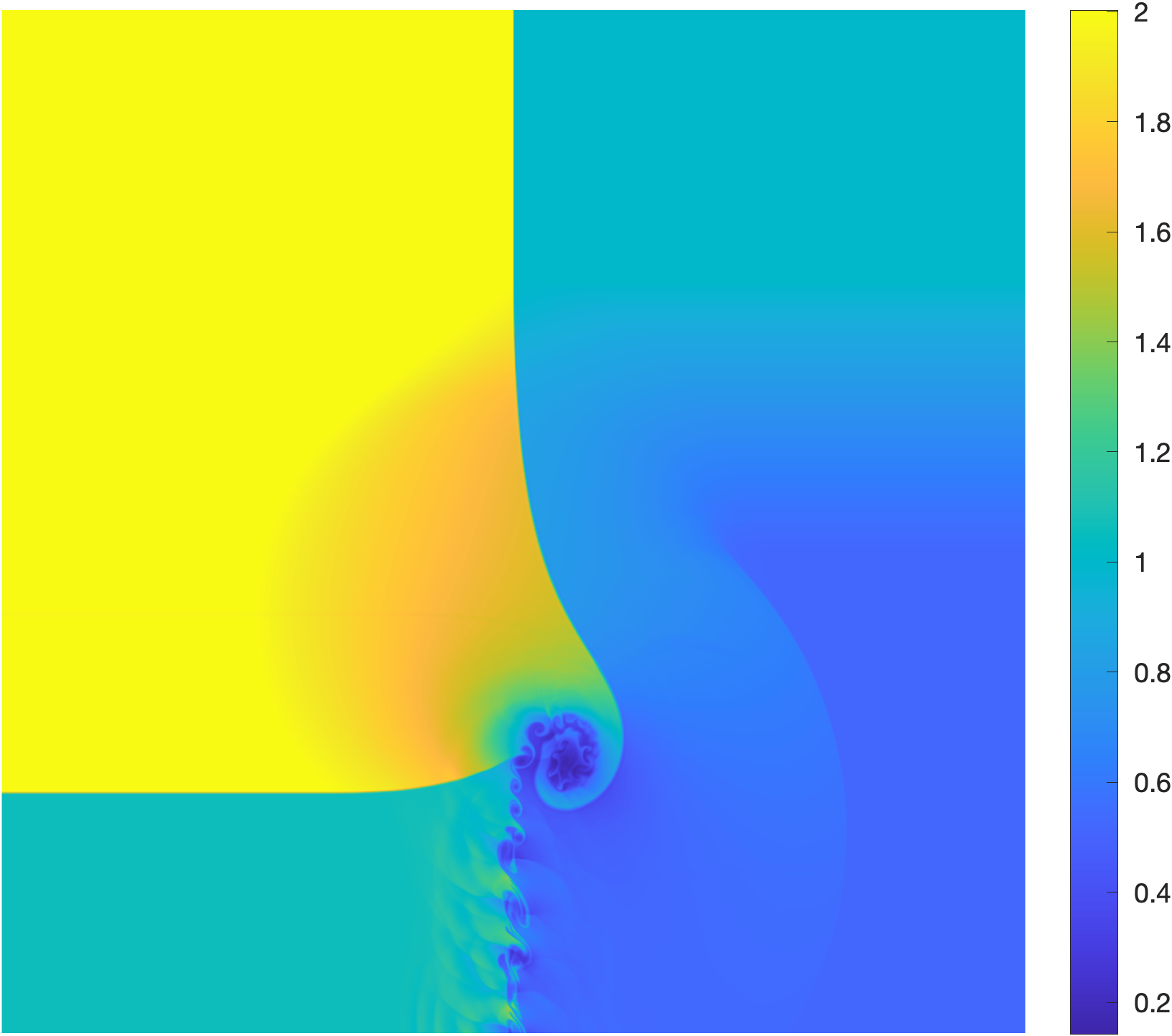}
	\caption{\label{fig:C17} Numerical results for
		Configuration 17 with initial perturbation for $x=x_0$, $y<y_0$ on grids with $256^2$ (left), $512^2$ (middle), and
		$1024^2$ (right) grid cells.
		Plot of density at time $t=0.3$ using
		moving-grid approach (top) and EG2 (bottom).}
\end{figure}
While the Kelvin-Helmholtz
instability develops for both versions of the Active Flux method,
it already appears on coarser grids when the moving-grid approach is used.

Next, we consider the Gresho vortex for different Mach numbers,
a standard problem for assessing the low Mach number properties of
numerical methods.
The velocity field of this stationary solution is divergence-free,
and the Mach number is varied by changing the background pressure.

\begin{example}
	We approximate solutions of the Gresho problem, a stationary rotating vortex. The initial values are given by
	\begin{equation*}
		\mathbf{u}(x,y,0)=
		\mathbf{n}
		\begin{cases}
			5r, & 0 \le r < 0.2,\\
			2-5r, & 0.2 \le r < 0.4,\\
			0, & r \ge 0.4,
		\end{cases}
	\end{equation*}
	and
	\begin{equation*}
		p(x,y,0)=
		\begin{cases}
			p_0+\dfrac{25}{2}r^2,
			& 0 \le r < 0.2,\\[0.5em]
			p_0+4-4\ln(0.2)
			+\dfrac{25}{2}r^2
			-20r
			+4\ln(r),
			& 0.2 \le r < 0.4,\\[0.5em]
			p_0-2+4\ln(2),
			& r \ge 0.4,
		\end{cases}
	\end{equation*}
	with
	\begin{equation*}
		\rho(x,y,0)=1,
		\qquad
		p_0=\frac{1}{\gamma Ma^2}.
	\end{equation*}
	Furthermore,
	\begin{equation*}
		r=\sqrt{x^2+y^2},
		\qquad
		\mathbf{n}
		=	(-\sin\theta,\cos\theta)^T
	\end{equation*}
	where $\theta\in[0,2\pi)$ denotes the polar angle and
	$\mathbf{u}=(u,v)^T$.
	
	The computational domain is $[-0.5,0.5]^2$, and periodic boundary conditions are imposed in both the $x$- and $y$-directions.
\end{example}

\begin{figure}[htb] 
	\centering
	\includegraphics[trim=0.5cm 0.5cm 0.5cm 0.cm,clip,width=0.32\textwidth]{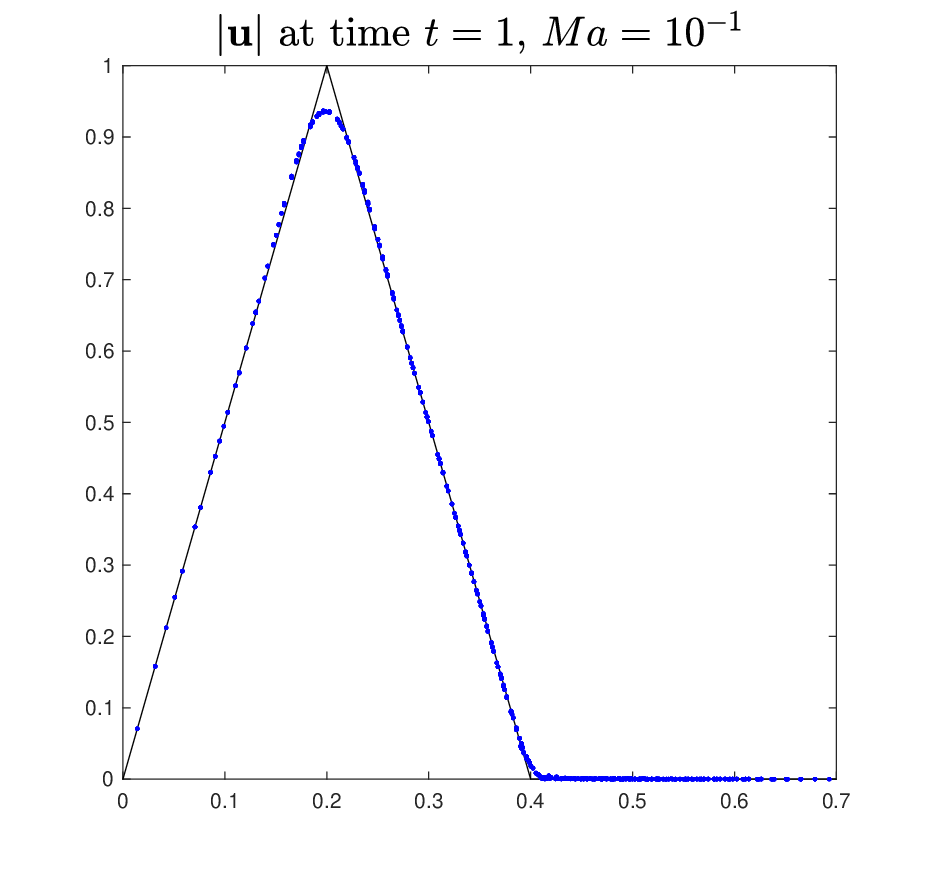}
	\hfill
	\includegraphics[trim=0.5cm 0.5cm 0.5cm 0.cm,clip,width=0.32\textwidth]{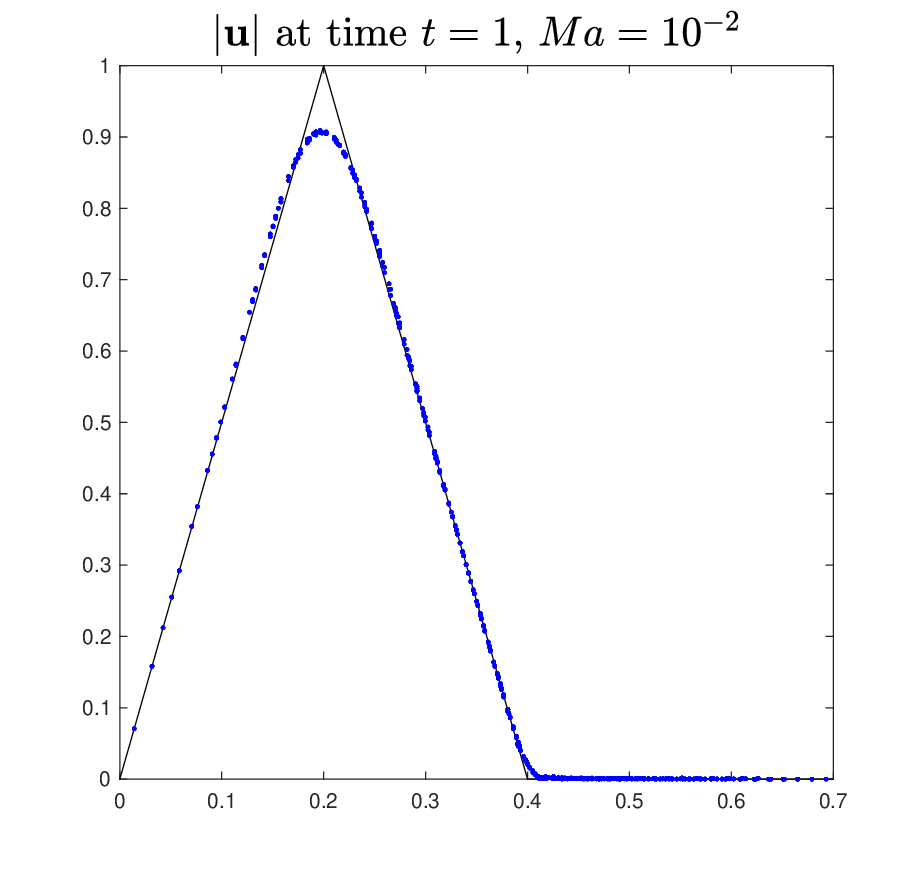}
	\hfill
	\includegraphics[trim=0.5cm 0.5cm 0.5cm
	0.cm,clip,width=0.32\textwidth]{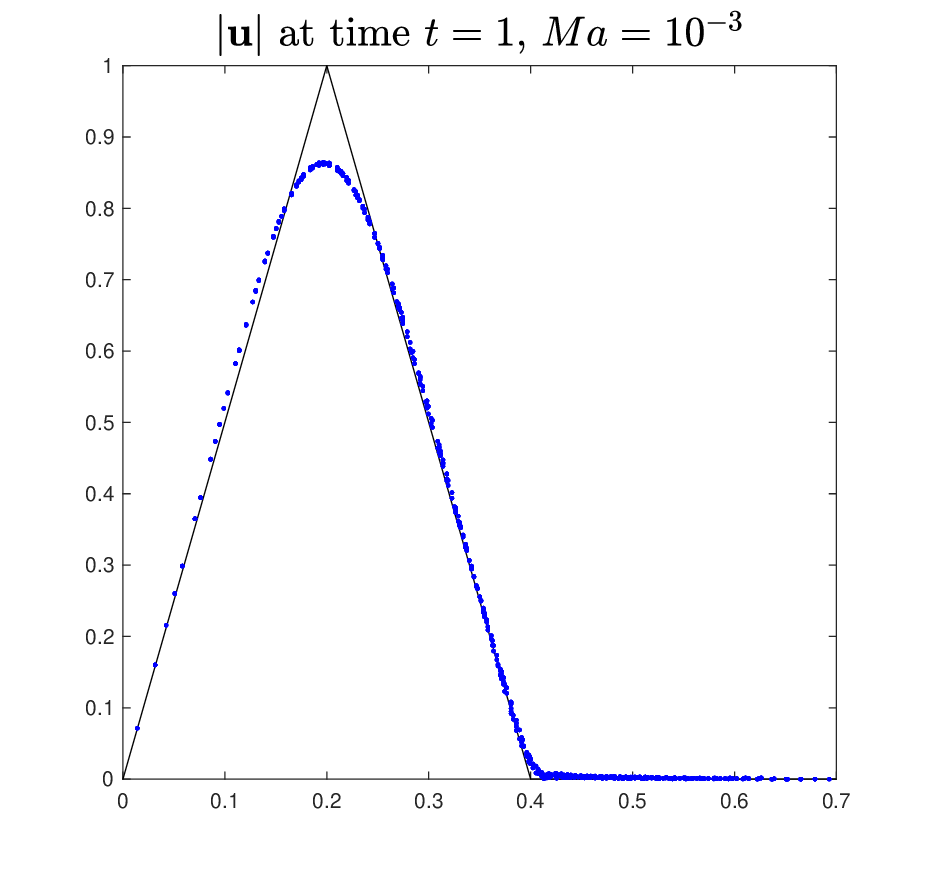}
	\\
	\includegraphics[trim=0.5cm 0.5cm 0.5cm 0.cm,clip,width=0.32\textwidth]{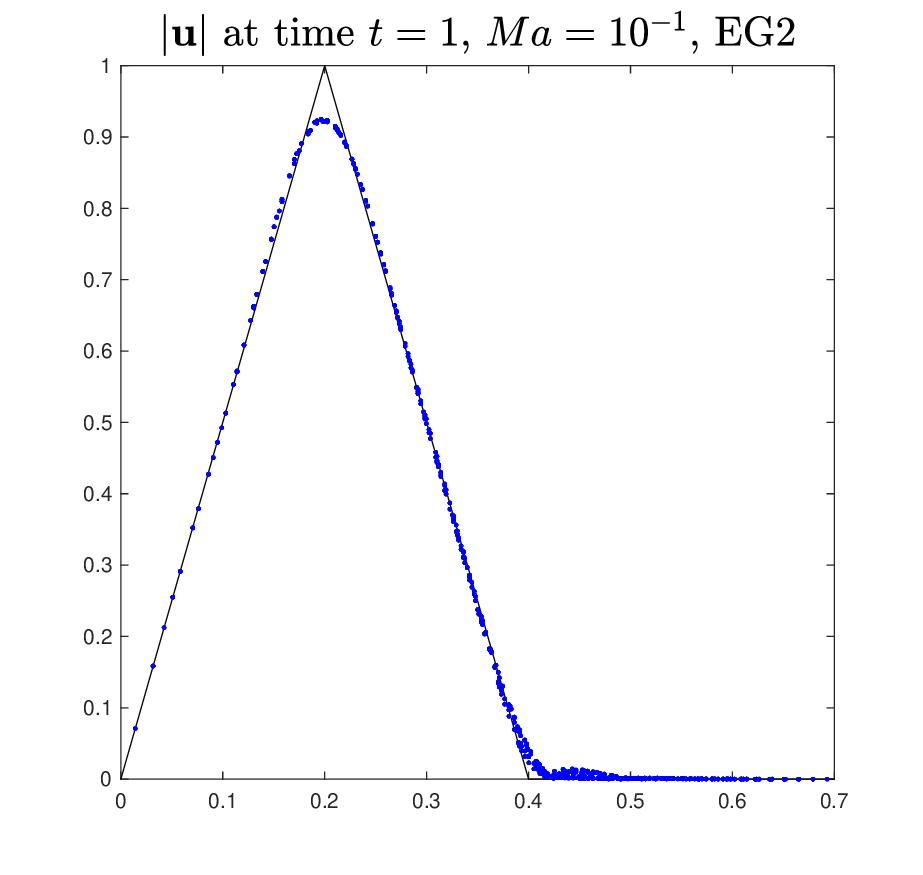}
	\hfill
	\includegraphics[trim=0.5cm 0.5cm 0.5cm 0.cm,clip,width=0.32\textwidth]{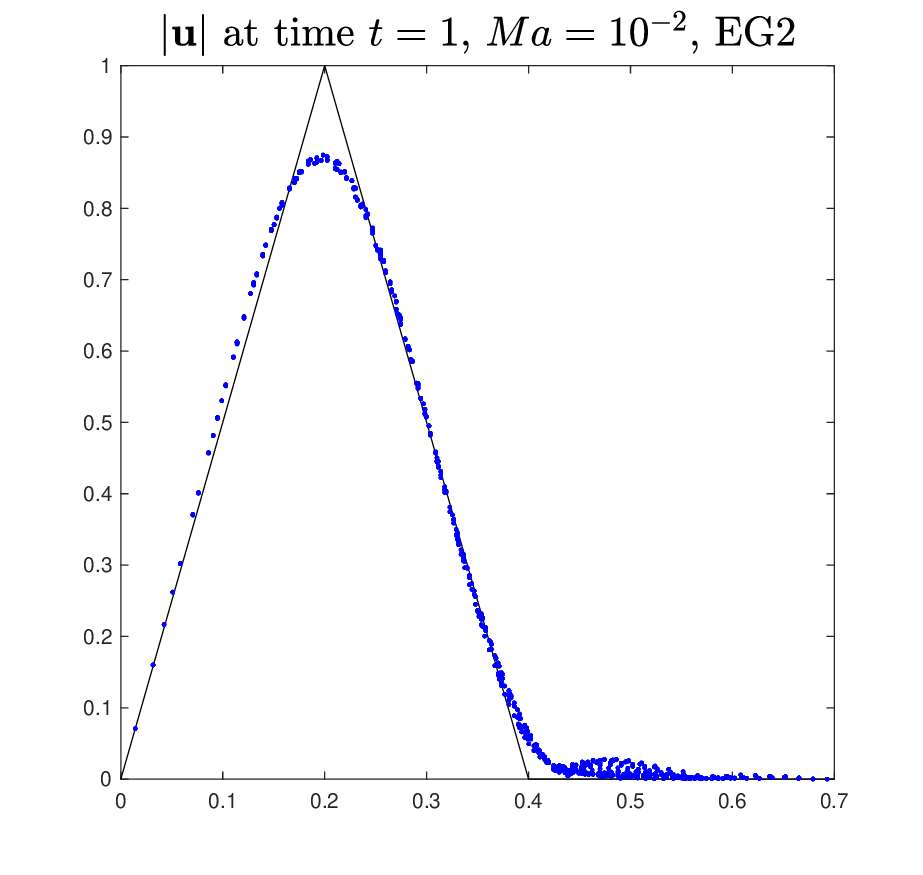}
	\hfill
	\includegraphics[trim=0.5cm 0.5cm 0.5cm 0.cm,clip,width=0.32\textwidth]{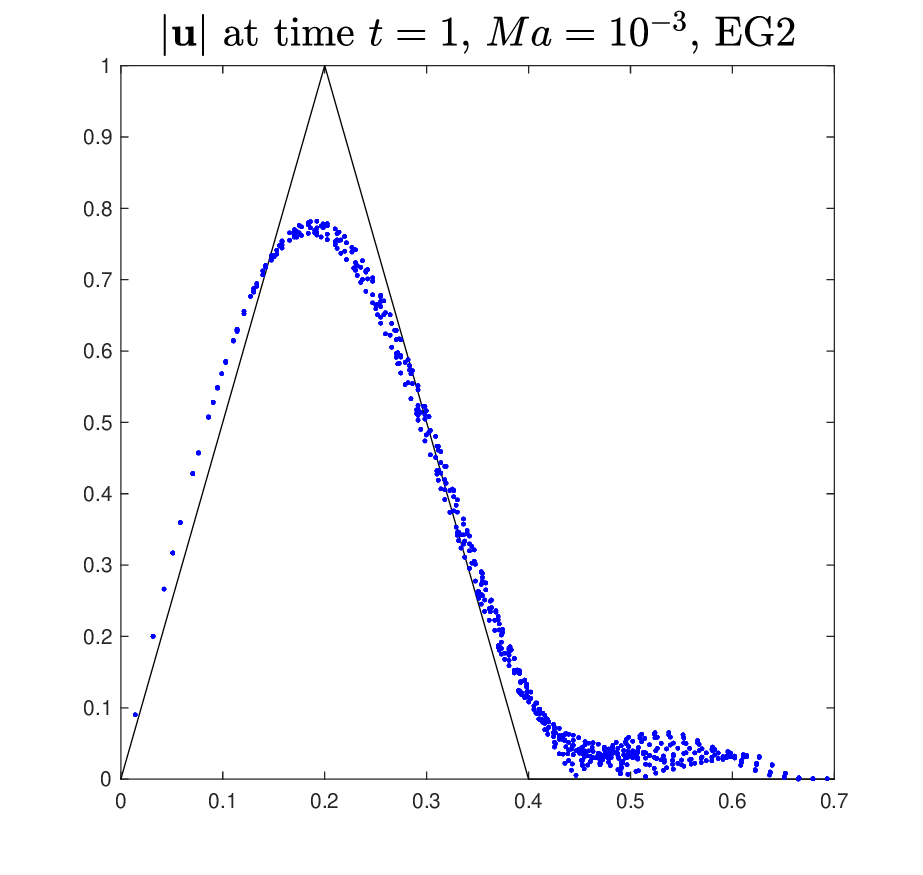} \\
	\caption{\label{fig:Gresho} Approximation of the Gresho vortex for different Mach
		numbers using the moving-grid approach (top) and the
		EG2 evolution operator (bottom).}
\end{figure}

In Figure \ref{fig:Gresho}, we show scatter plots of $\vert \mathbf{u}\vert$
at time $t = 1$ using different Mach numbers, i.e., $Ma = 10^{-1} $,
$Ma = 10^{-2} $ and $Ma = 10^{-3} $. We used a coarse grid with
$50^2$ cells. While both methods perform comparably for $Ma =
10^{-1}$, differences become apparent as the Mach number decreases.
In this regime, the method employing the exact evolution
operator preserves the vortex structure more accurately.

\begin{figure}[htb]
	\centering
	\begin{tabular}{cc}
		\includegraphics[width=0.45\textwidth]{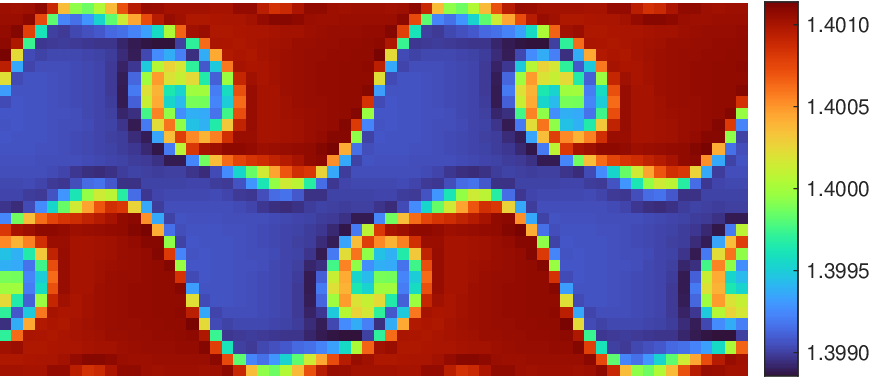} &
		\includegraphics[width=0.45\textwidth]{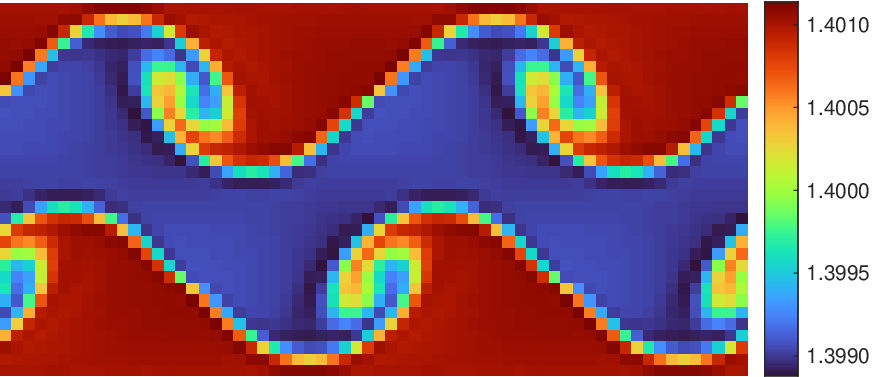}
		\\[0.15cm]
		\includegraphics[width=0.45\textwidth]{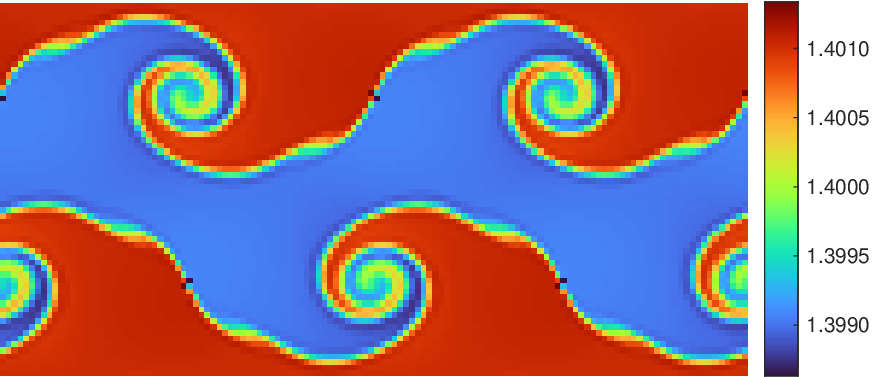} &
		\includegraphics[width=0.45\textwidth]{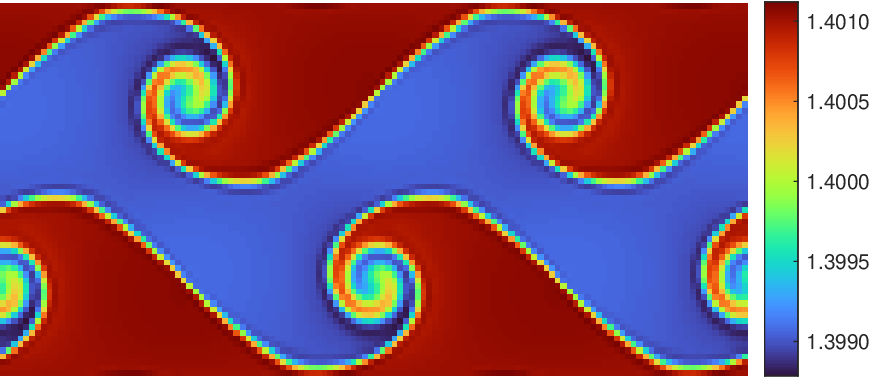} \\[0.15cm]
		\includegraphics[width=0.45\textwidth]{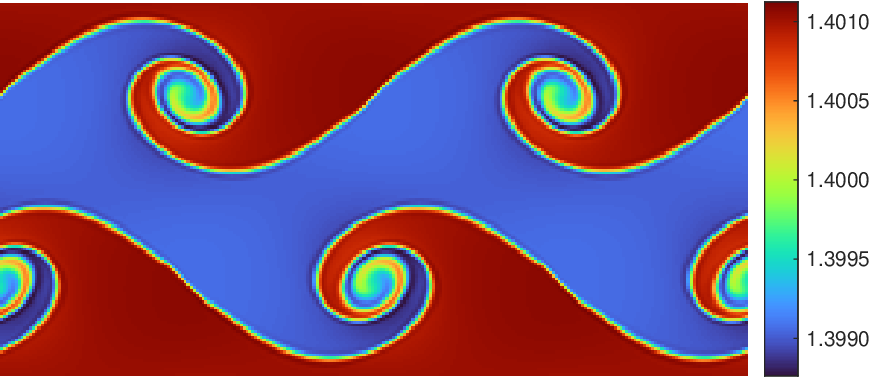} &
		\includegraphics[width=0.45\textwidth]{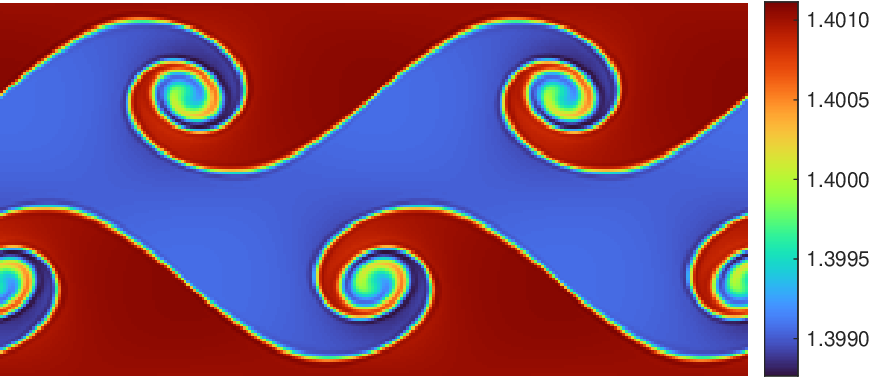}
	\end{tabular}
	\caption{Approximation of the Kelvin-Helmholtz instability on grids with
		$64\times 32$ cells (top), $128\times 64$ cells (middle), and
		$256\times128$ cells (bottom) using the moving
		grid approach (left) and EG2 (right).\label{fig:KH2}}
\end{figure}

The next example was proposed in \cite{article:LABHER2024} and also
considered in \cite{article:Barsukow2025,duraisamy2026}.
It describes a Kelvin-Helmholtz instability, which allows a
qualitative assessment of the numerical diffusion of a method.
Since the initial vorticity is smooth, the solution exhibits, for a
certain time, a well-defined evolution to which numerical methods
converge upon mesh refinement.
The resulting vortical structures provide a useful indicator of the
dissipative properties of a numerical method.
\begin{example}\label{ex:KH2}
	We consider the two-dimensional Euler equations with initial values of
	the form
	\begin{align*}
		\rho _0(x,y) = &\gamma + R(1-2\eta(y)),\\
		u _0(x,y) =& M(1-2\eta(y)),\\
		v_0 (x,y) = & \delta M \sin(2\pi x),\\		
		p_0 (x,y)  = & 1,
	\end{align*}
	with 
	\begin{equation*}
		\eta (y)= 
		\left\{ \begin{array}{ll}
			\frac{1}{2}(1+\sin (16 \pi (y+\frac{1}{4}))),  &-\frac{9}{32}\le y<-\frac{7}{32},\\
			1,  &-\frac{7}{32}\le y<\frac{7}{32},\\
			\frac{1}{2}(1+\sin (16 \pi (y+\frac{1}{4}))), &\frac{7}{32}\le y<\frac{9}{32},\\
			0,  &\text{else,}\\\end{array}\right.
	\end{equation*}
	and $R=10^{-3}$, $\delta = 0.1$ and $M = 0.01$, which is governing the Mach number of the flow. The domain is $[0,2]\times[-\frac{1}{2},\frac{1}{2}]$ with periodic boundary conditions in $x$- and $y$- directions. 
\end{example}
In Figure \ref{fig:KH2}, we show the results at $t=80$ calculated on
grids with different resolution.
All of these computations have been performed without limiting.
On the coarser grids, the approximations are visibly different, and the
solution obtained with the moving-grid approach appears to be
more deformed in some regions away from the vortices. Such deformations  
were also observed in coarse grid approximations 
shown in \cite{article:Barsukow2025,duraisamy2026}.
The Active Flux method with the EG2 evolution operator
provides an accurate approximation of the overall solution
structure while the vortices on the coarse grid are slightly less
rolled up.
It seems that achieving a more accurate approximation of the vortex
structure on a very coarse grid comes at the cost of a less
accurate approximation of the global solution structure. 
The solution structure of both methods compares
very well on the finest grid,
including the vortex approximation.

\section*{Conclusions}
We presented our latest developments on fully discrete, third-order
accurate Cartesian grid
Active Flux methods for the Euler equations of gas dynamics.
A crucial component of
these methods is the evolution operator, which is used to update the
point value degrees of freedom. Here we compare our recently proposed
Active Flux method, which uses an approximate evolution operator derived
using the method of bicharacteristics, with a new moving-grid approach
that provides an exact evolution operator for the linearized Euler
equations. 
The latter approach provides a more accurate approximation
of the stationary Gresho vortex in the low Mach number regime.
This fits well to previous results of Barsukow et al. \cite{article:Barsukow2020,article:BHKR2019}, which
showed that the Cartesian grid Active Flux method with the exact evolution
operator for acoustics is vorticity preserving. Solution structures
obtained with the moving-grid approach compare well with the solution
structures obtained by Barsukow's splitting approach with the exact
acoustic operator \cite{article:Barsukow2025}.

In other tests, we also observed a more pronounced
approximation of vortical structures when using the moving-grid
approach. We proposed a modification of a classical two-dimensional
Riemann problem where a stationary contact discontinuity was initially
slightly perturbed to trigger a Kelvin-Helmholtz instability. The
moving-grid approach exhibited a noticeably stronger growth of the
instability already on a coarser mesh than the Active Flux method
using the EG2 evolution operator.
In addition to this observation, both versions of the fully discrete Active Flux method
provided accurate approximations of
classical two-dimensional Riemann problems, even when relatively coarse
grids were used. These numerical experiments also suggest that the moving-grid approach requires somewhat stronger limiting than the Active Flux method based on the EG2 evolution operator. In particular, careful limiting of the local linearization states, together with the limiting of the point values, was found to be important for obtaining robust approximations near shocks. While the proposed limiting strategy in \cite{article:CHP2025} was effective for both methods, some problems involving shocks could be solved without limiting when the EG2 evolution operator was used, whereas the moving-grid approach required limiting.

Underresolved approximations of the
Kelvin-Helmholtz instability, as proposed in
\cite{article:LABHER2024}, indicated
that the Active Flux method with the EG2 evolution operator may provide a
more accurate overall solution structure when the solution consists
not only of isolated vortices. At the same time, the moving-grid
approach leads to a more pronounced approximation of the 
vortices, which on coarse grids required compensation in the form of a
less accurate approximation away from the vortices.

We believe that our studies can contribute to a better understanding of
underresolved approximations of turbulent solution structures.

\section*{End Matter}

\funding{This work was funded by the Deutsche Forschungsgemeinschaft (DFG, German Research Foundation) through project 525800857.}

\complia{On behalf of all authors, the corresponding author states that there is no conflict of interest.}

\bibliographystyle{camc}
\bibliography{sample}
\end{document}